\documentclass{birkmult}
\usepackage{amscd}
\usepackage{amsmath}
\usepackage{amsxtra}
\usepackage{amsfonts}
\usepackage{amssymb}

\newtheorem{theorem}{Theorem}[subsection]
\newtheorem{corollary}[theorem]{Corollary}
\newtheorem{lemma}[theorem]{Lemma}
\newtheorem{proposition}[theorem]{Proposition}
\theoremstyle{definition}
\newtheorem{definition}[theorem]{Definition}
\newtheorem{remark}[theorem]{Remark}

\newtheorem{example}[theorem]{Example}
\theoremstyle{remark}

\renewcommand{\theclaim}{\textup{\theclaim}}

\newtheorem*{acknowledgements}{Acknowledgements}

\newlength{\displayboxwidth}
\renewcommand{\theenumi}{\roman{enumi}}
\renewcommand{\labelenumi}{(\theenumi)}
\newcounter{saveenumi}
\numberwithin{equation}{section}
\def\openone

{\mathchoice

{\hbox{\upshape \small1\kern-3.3pt\normalsize1}}

{\hbox{\upshape \small1\kern-3.3pt\normalsize1}}

{\hbox{\upshape \tiny1\kern-2.3pt\SMALL1}}

{\hbox{\upshape \Tiny1\kern-2pt\tiny1}}}

\makeatletter

\newbox\ipbox

\newcommand{\ip}[2]{\left\langle\,#1\mid#2\,\right\rangle}

\newcommand{\diracb}[1]{\left\langle #1\mathrel{\mathchoice

{\setbox\ipbox=\hbox{$\displaystyle \left\langle\mathstrut
#1\right.$}

\vrule height\ht\ipbox width0.25pt depth\dp\ipbox}

{\setbox\ipbox=\hbox{$\textstyle \left\langle\mathstrut
#1\right.$}

\vrule height\ht\ipbox width0.25pt depth\dp\ipbox}

{\setbox\ipbox=\hbox{$\scriptstyle \left\langle\mathstrut
#1\right.$}

\vrule height\ht\ipbox width0.25pt depth\dp\ipbox}

{\setbox\ipbox=\hbox{$\scriptscriptstyle \left\langle\mathstrut
#1\right.$}

\vrule height\ht\ipbox width0.25pt depth\dp\ipbox}

}\right. }

\newcommand{\dirack}[1]{\left. \mathrel{\mathchoice

{\setbox\ipbox=\hbox{$\displaystyle \left.\mathstrut
#1\right\rangle$}

\vrule height\ht\ipbox width0.25pt depth\dp\ipbox}

{\setbox\ipbox=\hbox{$\textstyle \left.\mathstrut
#1\right\rangle$}

\vrule height\ht\ipbox width0.25pt depth\dp\ipbox}

{\setbox\ipbox=\hbox{$\scriptstyle \left.\mathstrut
#1\right\rangle$}

\vrule height\ht\ipbox width0.25pt depth\dp\ipbox}

{\setbox\ipbox=\hbox{$\scriptscriptstyle \left.\mathstrut
#1\right\rangle$}

\vrule height\ht\ipbox width0.25pt depth\dp\ipbox}

} #1\right\rangle}

\newcommand{\ltwor}{L^{2}\left(\mathbb{R}\right)}

\newcommand{\cj}[1]{\overline{#1}}

\newcommand{\bz}{\mathbb{Z}}
\newcommand{\Br}{\mathbb{R}}
\newcommand{\bc}{\mathbb{C}}
\newcommand{\bt}{\mathbb{T}}
\newcommand{\bn}{\mathbb{N}}

\newcommand{\xir}{X_\infty(r)}
\usepackage{graphicx}

\begin{document}
\title[Multiscale Theory and Wavelets in Nonlinear Dynamics]
{Methods from Multiscale Theory and Wavelets Applied to Nonlinear Dynamics}
\author[Dutkay]{Dorin Ervin Dutkay}

\address{%
Department of Mathematics\\
Rutgers University\\
New Brunswick, NJ 08901\\
U.S.A.}

\email{ddutkay@math.rutgers.edu}

\author[Jorgensen]{Palle E. T. Jorgensen}
\address{Department of Mathematics\\
The University of Iowa\\
Iowa City, IA 52242\\
U.S.A.}

\email{jorgen@math.uiowa.edu}

\thanks{Research supported in part by the National Science Foundation}
\begin{abstract}
We show how fundamental ideas from signal processing,
multiscale theory and wavelets may be applied to nonlinear dynamics.
\par
The problems from dynamics include iterated function systems (IFS),
dynamical systems based on substitution such as the discrete systems built
on rational functions of one complex variable and the corresponding Julia
sets, and state spaces of subshifts in symbolic dynamics. Our paper
serves to motivate and survey our recent results in this general area.
Hence we leave out some proofs, but instead add a number of intuitive
ideas which we hope will make the subject more accessible to researchers
in operator theory and systems theory.
\end{abstract}
\maketitle
\section{\label{S1}Introduction}
\par
In the past twenty years there has been a great amount of interest
in the theory of wavelets, motivated by applications to various
fields such as signal processing, data compression, tomography,
and subdivision algorithms for graphics (Our latest check on the
word ``wavelet'' in Google turned up over one million and a half
results, 1,590,000 to be exact). It is enough here to mention two
outstanding successes of the theory: JPEG 2000, the new standard
in image compression, and the WSQ (wavelet scalar quantization)
method which is used now by the FBI to store its fingerprint
database. As a mathematical subject, wavelet theory has found
points of interaction with functional and harmonic analysis,
operator theory, ergodic theory and probability, numerical
analysis and differential equations. With the explosion of
information due to the expansion of the internet, there came a need
for faster algorithms and better compression rates. These have
motivated the research for new examples of wavelets and new
wavelet theories.

\par
 Recent developments in wavelet analysis have brought together ideas
from engineering, from computational mathematics, as well as fundamentals
from representation theory. This paper has two aims: One to stress the
interconnections, as opposed to one aspect of this in isolation; and
secondly to show that the fundamental Hilbert space ideas from the linear
theory in fact adapt to a quite wide class of nonlinear problems. This class
includes random-walk models based on martingales. As we show, the theory
is flexible enough to allow the adaptation of pyramid algorithms to
computations in the dynamics of substitution theory as it is used in
discrete dynamics; e.g., the complex dynamics of Julia sets, and of
substitution dynamics.

\par

Our subject draws on ideas from a variety of directions. Of these
directions, we single out quadrature-mirror filters from signal/image
processing. High-pass/low-pass signal processing algorithms have now been
adopted by pure mathematicians, although they historically first were
intended for speech signals, see \cite{MR2004g:42041}. Perhaps unexpectedly, essentially
the same quadrature relations were rediscovered in operator-algebra
theory, and they are now used in relatively painless constructions of
varieties of wavelet bases. The connection to signal processing is rarely
stressed in the math literature. Yet, the flow of ideas between signal
processing and wavelet math is a success story that deserves to be told.
Thus, mathematicians have borrowed from engineers; and the engineers may
be happy to know that what they do is used in mathematics.
\par

 Our presentation serves simultaneously to motivate and to survey a
number of recent results in this general area. Hence we leave out some
proofs, but instead we add a number of intuitive ideas which we hope will
make the subject more accessible to researchers in operator theory and
systems theory. Our theorems with full proofs will appear elsewhere. An
independent aim of our paper is to point out several new and open
problems in nonlinear dynamics which could likely be attacked with the
general multiscale methods that are the focus of our paper.
\par

In Section \ref{S2} below, we present background material from signal
processing and from wavelets in a form which we hope will be accessible
to operator theorists, and to researchers in systems theory. This is
followed
in Section \ref{Exa}
by a presentation of some motivating examples from nonlinear
dynamics. They are presented such as to allow the construction of
appropriate Hilbert spaces which encode the multiscale structure.
Starting with a state space $X$ from dynamics, our main tool in the
construction of a multiscale Hilbert space $H(X)$ is the theory of random
walk and martingales.
The main part of our paper, Sections \ref{ProPro} and \ref{other},
serves to present our
recent new and joint results.

\section{\label{S2}Connection to signal processing and wavelets}
\renewcommand{\thetheorem}{\thesection.\arabic{theorem}}
We will use the term ``filter'' in the sense of signal processing.
In the simplest case, time is discrete, and a signal is just a
numerical sequence. We may be acting on a space of signals
(sequences) using the operation of Cauchy product; and the
operations of down-sampling and up-sampling. This viewpoint is
standard in the engineering literature, and is reviewed in \cite{BrJo02b} (see also \cite{BaVi02} and \cite{MR2003d:46102})
for the benefit of mathematicians.
\par
A numerical sequence $(a_k)$ represents a filter, but it is
convenient, at the same time, also to work with the frequency-%
response function. By this we mean simply the Fourier series (see (\ref{eqlowpass}) below)
corresponding to the sequence $(a_k)$. This Fourier series is
called the filter function, and in one variable we view it as a
function on the one-torus. (This is using the usual identification
of periodic functions on the line with functions on the torus.)
The advantage of this dual approach is that Cauchy product of
sequences then becomes pointwise product of functions.
\par
        We will have occasion to also work with several dimensions $d$, and
then our filter function represents a function on the $d$-torus. While signal
processing algorithms are old, their use in wavelets is of more recent
vintage. From the theory of wavelets, we learn that filters are used in the
first step of a wavelet construction, the so-called multiresolution
approach. In this approach, the problem is to construct a function on $\Br$ or
on $\Br^d$ which satisfies a certain renormalization identity, called the
scaling equation, and we recall this equation in two versions (\ref{eqscaling0}) and
(\ref{eqscaling}) below. The numbers $(a_k)$ entering into the scaling equation turn out to
be the very same numbers the signal processing engineers discovered as
{\it quadrature-mirror filters}.
\par
A class of wavelet bases are determined by filter functions. We
establish a duality structure for the category of filter functions and
classes of representations of a certain $C^*$-algebra. In the process, we
find new representations, and we prove a correspondence which serves as a
computational device. Computations are done in the sequence space $\ell^2$,
while wavelet functions are in the Hilbert space $L^2(\mathbb{R})$, or some other
space of functions on a continuum. The trick is that we choose a subspace
of $L^2(\mathbb{R})$ in which computations are done. The choice of subspace is
dictated by practical issues such as resolution in an image, or
refinement of frequency bands.

\par

 We consider non-abelian algebras containing canonical maximal abelian subalgebras,
i.e., $C(X)$ where $X$ is the Gelfand space, and the representations define
measures $\mu$ on $X$. Moreover, in the examples we study, it turns out that $X$ is an affine iterated function system (IFS), of course depending on the representation. In the
standard wavelet case, $X$ may be taken to be the unit interval. Following
Wickerhauser et al.\ \cite{CMW95}, these measures $\mu$ are used in the characterization and
analysis of wavelet packets.

\par
Orthogonal wavelets, or wavelet frames, for $L^{2}\left(  \mathbb{R}^{d}\right)  $
are associated with quadrature-mirror filters (QMF), a set of complex numbers
which relate the dyadic scaling of functions on $\mathbb{R}^{d}$ to the
$\mathbb{Z}^{d}$-translates. In the paper
\cite{MR2002h:46092}, 
we show that
generically, the data in the QMF-systems of wavelets are minimal, in the sense
that it cannot be nontrivially reduced. The minimality property is given a
geometric formulation in $\ell^{2}\left(  \mathbb{Z}^{d}\right)  $; and it is then
shown that minimality corresponds to irreducibility of a wavelet
representation of the Cuntz algebra $\mathcal{O}_{N}$. Our result is that this
family of representations of $\mathcal{O}_{N}$ on $\ell^{2}\left(
\mathbb{Z}^{d}\right)  $ is irreducible for a generic set of values of the
parameters which label the wavelet representations. Since MRA-wavelets
correspond to representations of the Cuntz algebras $\mathcal{O}_{N}$, we then
get, as a bonus, results about these representations.

\begin{definition} (Ruelle's operators.)
Let $(X, \mu)$ be a finite
Borel measure space, and let $r \colon X \rightarrow X$ be a
finite-to-one mapping of $X$ onto itself. (The measure $\mu$ will be
introduced in Section \ref{ProPro} below, and it will be designed so as to
have a certain strong invariance property which generalizes the
familiar property of Haar measure in the context of compact
groups.)
Let $V \colon X \rightarrow [0, \infty)$ be a given
measurable function. We then define an associated operator $R =
R_V$, depending on both $V$ and the endomorphism $r$, by the
following formula
\begin{equation}
R_Vf(x)=\frac{1}{\#\, r^{-1}(x)}\sum_{r(y)=x}V(y)f(y),\qquad f\in
L^1(X,\mu).
\end{equation}
Each of the operators $R_V$ will be referred to as a Ruelle operator, or a
transition operator; and each $R_V$ clearly maps positive functions to
positive functions. (When we say ``positive'' we do not mean
``strictly positive'', but merely ``non-negative''.)
\end{definition}
\par

We refer to our published papers/monographs \cite{Jor01a,Dut02,DuJo03} about the spectral picture of a positive transition operator $R$
(also called the Ruelle operator, or the Perron--Frobenius--Ruelle operator).
Continuing
\cite{BrJo99a}, 
it is shown in
\cite{Jor01a} 
that a general family of Hilbert-space constructions, which includes the
context of wavelets, may be understood from the so-called Perron--Frobenius
eigenspace of $R$. This eigenspace is studied further in the book
\cite{BrJo02b} 
by Jorgensen and O.~Bratteli;
see also
\cite{MR2004g:42041}. 
It is shown in
\cite{Jor01a} 
that a \emph{scaling identity}
(alias \emph{refinement equation})
\begin{equation}\label{eqscaling0}
\varphi\left(  x\right)  =\sqrt{N}\sum_{k\in\mathbb{Z}}a_{k}\varphi\left(
Nx-k\right)
\end{equation}
may be viewed generally and abstractly.
Variations of 
(\ref{eqscaling0}) open up for a variety of
new applications which we outline briefly below.
\renewcommand{\thetheorem}{\thesubsection.\arabic{theorem}}
\section{\label{Exa}Motivating examples, nonlinearity}

In this paper we outline how basic operator-theoretic notions from wavelet theory (such as the succesful multiresolution construction) may be adapted to certain state spaces in nonlinear dynamics. We will survey and extend our recent papers on multiresolution analysis of state spaces in symbolic shift dynamics $X(A)$, and on a class of substitution
systems $X(r)$ which includes the Julia sets; see, e.g., our papers
\cite{DuJo04c,DuJo04a,DuJo04b,DuJo03}. Our analysis of these systems $X$ starts with consideration of
Hilbert spaces of functions on $X$. But already this point demands new tools.
So the first step in our work amounts to building the right Hilbert spaces.
That part relies on tools from both operator algebras, and from probability
theory (e.g., path-space measures and martingales).

     First we must identify the appropriate measures on $X$, and secondly, we
must build Hilbert spaces on associated extension systems $X_\infty$, called
generalized solenoids.

     The appropriate measures $\mu$ on $X$ will be constructed using David
Ruelle's thermodynamical formalism \cite{Rue78a}: We will select our distinguished measures
$\mu$ on $X$ to have minimal free energy relative to a certain function $W$ on $X$.
The relationship between the measures $\mu$ on $X$ and a class of induced
measures on the extensions $X_\infty$ is based on a random-walk model which
we developed in \cite{DuJo04b}, such that the relevant measures on $X_\infty$ are
constructed as path-space measures on paths which originate on $X$. The
transition on paths is governed in turn by a given and prescribed function $W$
on $X$.

     In the special case of traditional wavelets, $W$ is the absolute
square of a so-called low-pass wavelet filter. In these special cases, the
wavelet-filter functions represent a certain system of functions on the
circle $\mathbb{T}=\mathbb{R}/\mathbb{Z}$  which we outline below; see also our paper \cite{DuJo03}.  Even our
analysis in \cite{DuJo03} includes wavelet bases on affine Cantor sets in $\Br^d$.
Specializing our Julia-set cases $X = X(r)$ to the standard wavelet systems,
we get $X = \mathbb{T}$, and the familiar approach to wavelet systems becomes a special
case.

\medskip

\par
According to its original definition, a {\it wavelet} is a function $\psi\in L^2(\mathbb{R})$ that generates an orthonormal basis under the action of two unitary operators: the dilation and the translation. In particular,
$$
\{\,U^jT^k\psi\mid j,k\in\mathbb{Z}\,\}
$$
must be an orthonormal basis for $L^2(\mathbb{R})$, where
\begin{equation}
Uf(x)=\frac{1}{\sqrt{2}}f\left(\frac{x}{2}\right),\,\quad Tf(x)=f(x-1),\qquad f\in L^2(\mathbb{R}),x\in\mathbb{R}.
\label{UT}
\end{equation}

One of the effective ways of getting concrete wavelets is to first look for
a \emph{multiresolution analysis} (\emph{MRA\/}), that is to first identify a suitable telescoping nest of
subspaces $(\mathcal{V}_n)_{n\in\mathbb{Z}}$ of $L^2(\mathbb{R})$ that have trivial intersection and dense union, $\mathcal{V}_n=U\mathcal{V}_{n-1}$, and $\mathcal{V}_0$ contains a {\it scaling function} $\varphi$ such that the translates of $\varphi$ form an orthonormal basis for $\mathcal{V}_0$.

The scaling function will
necessarily satisfy an equation of the form
\begin{equation}\label{eqscaling}
U\varphi=\sum_{k\in\mathbb{Z}}a_kT^k\varphi,\end{equation}
called the scaling equation (see also (\ref{eqscaling0})).

\par
To construct wavelets, one has to choose a {\it low-pass filter} \begin{equation}\label{eqlowpass}
m_0(z)=\sum_{k\in\mathbb{Z}}a_kz^k
\end{equation}
 and from this obtain the scaling function $\varphi$ as the fixed point of a certain cascade operator. The wavelet $\psi$ (mother function) is constructed from $\varphi$ (a father
function, or scaling function) with the aid of two closed subspaces $\mathcal{V}_0$, and
$\mathcal{W}_0$. The scaling function $\varphi$ is obtained as the solution to a scaling
equation (\ref{eqscaling}), and its integral translates generates $\mathcal{V}_0$ (the initial
resolution), and $\psi$ with its integral translates generating $\mathcal{W}_0:= \mathcal{V}_1
\ominus \mathcal{V}_0$, where $\mathcal{V}_1$ is the next finer resolution subspace containing
the initial resolution space $\mathcal{V}_0$, and obtained from a one-step refinement of
$\mathcal{V}_0$.
\par
Thus the main part of the construction involves a clever choice for the low-pass filter $m_0$ that gives rise to nice orthogonal wavelets.

\subsection{\label{MRAs}MRAs in geometry and operator theory}

Let $X$ be a compact\label{compactXa} Hausdorff space, and let $r\colon X\to X$ be a finite-to-one
continuous endomorphism, mapping $X$ onto itself. As an example, $r = r(z)$ could
be a rational mapping of the Riemann sphere, and $X$ could be the corresponding
Julia set; or $X$ could be the state space of a subshift associated to a $0$--$1$
matrix.

Due to work of Brolin\label{Brolina} \cite{MR33:2805} and Ruelle
\cite{Rue89}, it is known that for these examples, $X$ carries a
unique maximal entropy, or minimal free-energy measure $\mu$ also
called strongly $r$-invariant; see Lemma \ref{lemfme3_2} (\ref{lemfme3_2(2)}), and
equations (\ref{eqpr_4}) and (\ref{eq3_1})--(\ref{eq3_2}) below.
For each point $x \in X$, the measure $\mu$ distributes the ``mass'' equally
on the finite number of solutions $y$ to $r(y) = x$.

We show that this structure in fact
admits a rich family of wavelet bases. Now this will be on a
Hilbert space which corresponds to $L^2(\mathbb{R})$ in the familiar case of
multiresolution wavelets. These are the wavelets corresponding to scaling $x$ to
$Nx$, by a fixed integer $N$, $N \ge 2$. In that case, $X = \mathbb{T}$, the circle in $\mathbb{C}$, and $r
(z) = z^N$. So even in the ``classical'' case, there is a ``unitary dilation'' from
$L^2(X)$ to $L^2(\mathbb{R})$ in which the Haar measure on $\mathbb{T}$ ``turns into'' the Lebesgue
measure on $\mathbb{R}$.

Our work so far, on examples, shows that this viewpoint holds promise for
understanding the harmonic analysis of Julia sets, and of other iteration
systems. In these cases, the analogue of $L^2(\mathbb{R})$ involves quasi-invariant
measures $\mu_\infty$ on a space $X_\infty$, built from $X$ in a way that follows
closely the analogue of passing from $\mathbb{T}$ to $\mathbb{R}$ in wavelet theory. But $\mu_\infty$
is now a path-space measure. The translations by the integers $\mathbb{Z}$ on $L^2(\mathbb{R})$ are
multiplications in the Fourier dual. For the Julia examples, this corresponds to
the coordinate multiplication, i.e., multiplication by $z$ on $L^2(X,\mu)$, a normal
operator. We get a corresponding covariant system on $X_\infty$, where
multiplication by $f(z)$ is unitarily equivalent to multiplication by $f(r(z))$.
But this is now on the Hilbert space $L^2(X_\infty$), and defined from the 
path-space measure $\mu_\infty$.
Hence all the issues that are addressed since the mid-1980's for $L^2(\mathbb{R})$-%
wavelets have analogues in the wider context, and they appear to reveal
interesting spectral theory for Julia sets.

\par
\subsection{\label{ExaSpe}Spectrum and geometry: wavelets, tight frames, and Hilbert spaces on Julia sets}

\subsubsection{\label{SpectrumBackground}Background}

Problems in dynamics are in general irregular and chaotic, lacking
the internal
structure and homogeneity which is typically related to group actions. Lacking
are the structures that lie at the root of harmonic analysis, i.e., the torus
or $\mathbb{R}^d$.

  Chaotic attractors can be viewed simply as sets, and they are often lacking
the structure of manifolds, or groups, or even homogeneous manifolds.
Therefore, the ``natural'' generalization, and question to ask, then depends
on the point of view, or on the application at hand.

  Hence, many applied problems are typically not confined to groups; examples
are turbulence, iterative dynamical systems, and irregular sampling. But there
are many more. Yet harmonic analysis and more traditional MRA theory begins
with $\mathbb{R}^d$, or with some such group context.

  The geometries arising in wider contexts of applied mathematics might be attractors
in discrete dynamics, in substitution tiling problems, or in complex
substitution schemes, of which the Julia sets are the simplest.
These geometries are not tied to groups at all. And yet they are very
algorithmic, and they invite spectral-theoretic computations.

  Julia\label{Juliaa} sets are prototypical examples of chaotic attractors for iterative discrete
dynamical systems; i.e., for iterated substitution of a fixed rational function
$r(z)$, for $z$ in the Riemann sphere.
So the Julia sets serve beautifully as a testing ground for discrete
algorithms, and for analysis on attractors.
In our papers \cite{DuJo04b}, \cite{DuJo04c}
we show that multiscale/MRA analysis adapts itself well to discrete iterative
systems, such as Julia sets, and state space for subshifts in symbolic
dynamics.
But so far there is very little computational harmonic analysis outside the
more traditional context of $\mathbb{R}^d$.
\begin{definition}
Let\label{juliasetdefa}
$S$ be the Riemann sphere, and let $r\colon S \to S$ be a rational map of degree greater than one. Let $r^n$  be the $n$'th iteration of $r$,
i.e., the $n$'th iterated substitution of $r$ into itself. The Fatou set $F(r)$ of $r$ is the largest open set $U$ in $S$ such that the sequence $r^n$,
restricted to $U$, is a normal family (in the sense of Montel). The \emph{Julia set} $X(r)$ is the complement of $F(r)$.
\end{definition}
   Moreover, $X(r)$ is known to behave (roughly) like an attractor for the discrete dynamical system $r^n$. But it is more complicated: the
Julia set $X(r)$ is the locus of expanding and chaotic behavior; e.g., $X(r)$ is equal to the closure of the set of repelling periodic points
for the dynamical system $r^n$ in $S$.

In addition, $X = X(r)$ is the minimal closed set $X$, such that $|X| > 2$,  and $X$ is invariant under the branches of $r^{-1}$, the inverses of $r$.

  While our prior research on Julia sets is
only in the initial stages, there is already some progress. We are especially
pleased that Hilbert spaces of martingales have offered the essential framework
for analysis on Julia sets. The reason is that martingale tools are more adapted to
irregular structures than are classical Fourier methods.

To initiate a harmonic analysis outside the classical context of the tori and of
$\mathbb{R}^d$, it is thus natural to begin with the Julia sets. Furthermore, the Julia
sets are already of active current research interest in geometry and dynamical
systems.

But so far, there are only sporadic attempts at a harmonic analysis for Julia
sets, let alone the more general geometries outside the classical context of
groups.
Our prior work should prepare us well for this new project.
Now, for Julia-set geometries, we must begin with
the Hilbert space. And even that demands new ideas and new tools.

   To accommodate wavelet solutions in the wider context, we built new Hilbert
spaces with the use of tools from martingale theory. This is tailored to
applications
we have in mind
in geometry and dynamical systems involving Julia sets,
fractals, subshifts, or even more general discrete dynamical systems.

    The construction of solutions requires a combination of tools that are
somewhat non-traditional in the subject.

A number of wavelet constructions in pure and applied mathematics
have a common operator-theoretic underpinning. It may be illustrated
with the following operator system:
$\mathcal{H}$ some Hilbert space; $U\colon\mathcal{H}\rightarrow\mathcal{H}$ a unitary operator; $V\colon\mathbb{T}\rightarrow\mathcal{U}\left(  \mathcal{H}\right)  $ a
unitary representation. Here $\mathbb{T}$ is the $1$-torus, $\mathbb{T}%
=\mathbb{R}/\mathbb{Z}$, and $\mathcal{U}\left(  \mathcal{H}\right)  $ is the
group of unitary operators in $\mathcal{H}$.

The operator system satisfies the identity%
\begin{equation}
V\left(  z\right)  ^{-1}UV\left(  z\right)  =zU,\qquad z\in\mathbb{T}.
\label{eqPJ820.1}%
\end{equation}
\begin{definition}\label{defhomog}
We say that $U$ is homogeneous of degree one with respect to the
scaling automorphisms defined by $\left\{ \,V\left( z\right) \mid
z\in\mathbb{T}\,\right\}  $ if (\ref{eqPJ820.1}) holds.
 In addition, $U$
must have the spectral type defined by Haar measure. We say that $U$ is homogeneous of degree one if it satisfies (\ref{eqPJ820.1}) for some representation $V(z) $ of $\bt$.

\end{definition}
In the case of the standard wavelets with scale $N$, $\mathcal{H}=L^{2}\left(
\mathbb{R}\right)  $,
where in this case the two operators $U$ and $T$ are the usual $N$-adic scaling
and integral translation operators; see 
(\ref{UT}) above giving $U$
for the case $N = 2$.
\par
If $\mathcal{V}_{0}\subset\mathcal{H}$ is a resolution subspace, then
$\mathcal{V}_{0}\subset U^{-1}\mathcal{V}_{0}$. Set $\mathcal{W}_{0}%
:=U^{-1}\mathcal{V}_{0}\ominus\mathcal{V}_{0}$, $Q_{0}:={}$the projection onto
$\mathcal{W}_{0}$, $Q_{k}=U^{-k}Q_{0}U^{k}$, $k\in\mathbb{Z}$, and%
\begin{equation}
V\left(  z\right)  =\sum_{k=-\infty}^{\infty}z^{k}Q_{k},\qquad z\in\mathbb{T}.
\label{eqPJ820.4}%
\end{equation}
Then it is easy to check that the pair $\left(  U,V\left(  z\right)  \right)
$ satisfies the commutation identity (\ref{eqPJ820.1}).

Let $\left(  \psi_{i}\right)  _{i\in I}$ be a Parseval frame in $\mathcal{W}%
_{0}$. Then%
\begin{equation}
\left\{  \,U^{k}T_{n}\psi_{i}\bigm|i\in I,\;k,n\in\mathbb{Z}\,\right\}
\label{eqPJ820.5}%
\end{equation}
is a Parseval frame for $\mathcal{H}$. (Recall, a Parseval frame is also
called a normalized tight frame.)

Turning the picture around, we\label{systemUVa} may start with a system $\left(  U,V\left(
z\right)  \right)  $ which satisfies (\ref{eqPJ820.1}), and then reconstruct
wavelet (or frame) bases in the form (\ref{eqPJ820.5}). To do this in the
abstract, we must first understand the multiplicity function calculated for
$\left\{  \,T_{n}\mid n\in\mathbb{Z}\,\right\}  $ when restricted to the two
subspaces $\mathcal{V}_{0}$ and $\mathcal{W}_{0}$. But such a multiplicity
function may be be defined for any representation of an abelian $C^{\ast}%
$-algebra acting on $\mathcal{H}$ which commutes with some abstract scaling
automorphism ($V\left(  z\right)  $).

This technique can be used for the Hilbert spaces associated to Julia
sets, to construct wavelet (and frame) bases in this context.

\subsubsection{\label{SpectrumTech}Wavelet filters in
nonlinear models}

These systems are studied both in the theory of symbolic
dynamics, and in $C^*$-algebra theory, i.e., for the Cuntz--Krieger algebras.

It is known\label{Brolinb}
\cite{MR33:2805,Rue78a}
that, for these examples, $X$ carries a unique maximal entropy, or
minimal free-energy measure $\mu$. The most general case when $(X,r)$ admits a
strongly $r$-invariant measure is not well understood.

  The intuitive\label{massa}
property of such measures $\mu$ is this: For each point $x$ in $X$, $\mu$
distributes the ``mass'' equally on the finite number of solutions $y$ to $r(y) = x$.
Then $\mu$ is obtained by an iteration of this procedure,
i.e., by considering successively the finite sets of solutions $y$ to $r^n(y)=x$,
for all $n = 1,2, \dots$; taking an average, and then limit. While the procedure
seems natural enough, the structure of the convex set $K(X,r)$ of all the $r$-%
invariant measures is worthy of further study. For the case when $X =\bt$ the
circle, $\mu$ is unique, and it is the Haar measure on $\bt$.

  The invariant measures are of interest in several areas: in representation
theory, in dynamics, in operator theory, and in $C^*$-algebra theory. Recent
work of Dutkay--Jorgensen
\cite{DuJo04a,DuJo04b,DuJo04c}
focuses on special convex subsets of $K(X,r)$ and their
extreme points. This work in turn is motivated by our desire to construct
wavelet bases on Julia sets.
\par
Our work on $K(X,r)$ extends what is know in traditional wavelet analysis.
While linear wavelet analysis is based on Lebesgue measure on $\Br^d$, and Haar
measure on the $d$-torus $\bt^d$, the nonlinear theory is different, and it
begins with results from geometric measure theory: We rely on results from \cite{MR33:2805}, \cite{MR85m:58110a}, \cite{MR85m:58110b}, and \cite{MR2003772}. And our adaptation of spectral theory in turn is motivated by work by Baggett et al.\ in \cite{BJMP04,BJMP03b}.

Related work reported in \cite{MR2003772} further ties in with exciting advances on graph
$C^*$-algebras pioneered by Paul Muhly, Iain Raeburn and many more in the $C^*$-%
algebra community; see \cite{BJMP04}, \cite{BaVi02}, \cite{Bis95},
\cite{BJKR01}, \cite{BJKR02}, \cite{Con94}, \cite{MR2004i:46109}
\cite{Oka03}, \cite{Rad91}, \cite{Rie81}.

     We now prepare the tools in outline, function systems and self-similar
operators which will be used in Subsection \ref{MRonJ} in the construction of
specific Hilbert bases in the context of Julia sets.

Specifically, we outline a construction of a multiresolution/%
wavelet analysis for the Julia set $X$ of a given rational
function $r(z)$ of one complex variable, i.e.,
$X=\operatorname{Julia}(r)$. Such an analysis could likely be
accomplished with some kind of wavelet representation induced by a
normal operator $T$ which is unitarily equivalent to a function of
itself. Specifically, $T$ is unitarily equivalent to $r(T)$, i.e.,
\begin{equation}
UT=r(T)U\text{,\qquad for some unitary operator }U\text{ in a Hilbert space
}\mathcal{H}. \label{eq1}%
\end{equation}
Even the existence of these representations would seem to be new and
significant in operator theory.

There are consequences and applications of such a construction:
First we will get the existence of a finite set of some special
generating functions $m_{i}$ on the Julia set $X$ that may be
useful for other purposes, and which are rather difficult to
construct by alternative means. The relations we require for a
system $m_{0},\dots,m_{N-1}$ of functions on $X(r)$
are as follows:%
\begin{equation}
\frac{1}{N}\sum_{\substack{y\in X(r)\\r(y)=x}}\overline{m_{i}(y)}%
\,m_{j}(y)h(y)=\delta_{i,j}h(x)\text{\qquad for }\mathrm{a.e.}\,x\in X(r),
\label{eq2}%
\end{equation}
where $h$ is a Perron--Frobenius eigenfunction for a Ruelle
operator $R$ defined on $L^{\infty}(X(r))$. Secondly, it follows
that for each Julia set $X$, there is an infinite-dimensional
group which acts transitively on these systems of functions. The
generating functions on $X$ we are constructing are analogous to
the more familiar functions on the circle $\mathbb{T}$ which
define systems of filters in wavelet analysis. In fact, the
familiar construction of a wavelet basis in
$\mathcal{H}=L^{2}(\mathbb{R})$ is a special case of our analysis.

In the standard wavelet case, the rational function $r$ is just a monomial,
i.e., $r(z)=z^{N}$ where $N$ is the number of frequency bands, and $X$ is the
circle $\mathbb{T}$ in the complex plane. The simplest choice of the functions
$m_{0},\dots,m_{N-1}$ and $h$ in this case is $m_{k}(z)=z^{k}$, $z\in
\mathbb{T}$, $0\leq k\leq N-1$, and $h(z)\equiv1$, $z\in\mathbb{T}$. This
represents previous work by many researchers, and also joint work between
Jorgensen and Bratteli \cite{BrJo02b}.

Many applied problems are typically not confined to groups; examples
are turbulence, iterative dynamical systems, and irregular sampling. But there
are many more.

 Even though attractors such as Julia sets are highly nonlinear, there are
adaptations of the geometric tools from the analysis on $\mathbb{R}$ to
the nonlinear setting.

The adaptation of traditional wavelet tools to nonlinear settings begins
with an essential step: the construction of an appropriate Hilbert space;
hence the martingales.

Jorgensen, Bratteli, and Dutkay have developed a
representation-theoretic duality for wavelet filters and the associated
wavelets, or tight wavelet frames. It is based on representations of the Cuntz
algebra $\mathcal{O}_N$ (where $N$ is the scale number); see \cite{BrJo02b}, \cite{Dut04d}, \cite{DuJo04b}, \cite{MR2002h:46092}. As a by-product of this approach
they get infinite families of inequivalent representations of $\mathcal{O}_N$ which are of
independent interest in operator-algebra theory.


The mathematical community has long known of operator-theoretic power in Fourier-analytic issues,
e.g., Rieffel's incompleteness theorem \cite{Rie81}
for Gabor systems
violating the Nyquist condition, cf.\ \cite{DLL96}.
We now feel that we can address a number of significant problems combining operator-theoretic and harmonic analysis methods, as well as broadening and deepening the
mathematical underpinnings of our subject.
\subsection{\label{Multi}Multiresolution analysis \textup{(}MRA\/\textup{)}}

  One reason various geometric notions of multiresolution in Hilbert
space (starting with S. Mallat) have proved  useful in
computational mathematics is that these resolutions are modeled on the
fundamental concept of dyadic representation of the real numbers; or more
generally on some variant of the classical positional representation of real
numbers. In a dyadic representation of real numbers, the shift corresponds to
multiplication by $2^{-1}$.

     The analogue of that in Hilbert space is a unitary operator $U$ which
scales functions by $2^{-1}$. If the Hilbert space $\mathcal{H}$ is $L^2(\Br)$, then a
resolution is a closed subspace $\mathcal{V}_0$ in $\mathcal{H}$ which is invariant under $U$ with
the further property that the restriction of $U$ to $\mathcal{V}_0$ is a shift operator. Positive and negative powers of $U$ then scale $\mathcal{H}$ into a nested resolution system indexed by the integers; and this lets us solve algorithmic problems by simply following the standard rules from number theory.

     However, there is no reason at all that use of this philosophy should
be restricted to function spaces on the line, or on Euclidean space.
Other geometric structures from discrete dynamics
admit resolutions as well.

Most MRA,
Fourier multiresolution analysis (FMRA)
and generalized multiresolution analysis (GMRA)
wavelet constructions in pure and applied mathematics
have a common operator-theoretic underpinning. Consider the operator system
$\left(\mathcal{H},U,V\left(z\right)\right)$ from (\ref{eqPJ820.1}) above.

\begin{proposition}\label{defmulti}
Once the system\label{Otsa} $\left( \mathcal{H}, U,V\right)  $ is
given as in \textup{(\ref{eqPJ820.1})}, a scale of
closed subspaces $\mathcal{V}_{n}$ \textup{(}called resolution
subspaces\/\textup{)} may be derived from the
spectral subspaces of the representation
$\mathbb{T}\ni z\mapsto V\left(z\right)$, i.e., subspaces
$\left(\mathcal{V}_{n}\right)_{n\in\mathbb{Z}}$
such that $\mathcal{V}_{n}\subset\mathcal{V}_{n+1}$,
\[
\bigcap\mathcal{V}_{n}=\{0\},\qquad\bigcup\mathcal{V}_{n}\text{ is dense in }\mathcal{H},
\]
and $U\mathcal{V}_{n}\subset\mathcal{V}_{n-1}$. \end{proposition}
Conversely, if these spaces $\left(\mathcal{V}_{n}\right)$ are
given, then $V(z)$ defined by equation (\ref{eqPJ820.4}) can be
shown to satisfy (\ref{eqPJ820.1}) if
\[
Q_0={}\text{the orthogonal projection onto }U^{-1}\mathcal{V}_{0}\ominus\mathcal{V}_{0},
\]
and $Q_n=U^{-n}Q_0U^n$, $n\in\mathbb{Z}$.
\par As a result we note the following criterion.
\begin{proposition} A given unitary operator $U$ in a Hilbert space $\mathcal{H}$ is
part of some multiresolution system $(\mathcal{V}_n)$ for $\mathcal{H}$ if and only if $U$ is
homogeneous of degree one with respect to some representation $V(z)$ of the
one-torus $\bt$.
\end{proposition}
\par
In application the spaces $\mathcal{V}_{n}$ represent a grading of
the entire Hilbert space $\mathcal{H}$, and we say that $U$ scales
between the different levels. In many cases, this
operator-theoretic framework serves to represent structures which
are similar up to scale, structures that arise for example in the
study of \emph{fractional Brownian motion} (\emph{FBM}\/). See, e.g.,
\cite{Lin93}, which offers a Hilbert-space formulation of FBM based
on a white-noise stochastic integration.

Moreover, this can be done for the Hilbert space $\mathcal{H}=L^{2}\left(
X_{\infty},\mu_{\infty}\right)  $ associated with a Julia set $X$. As a
result, we get wavelet bases and Parseval wavelets in this context.

\par
This is relevant for the two mentioned classes of examples, Julia
sets $X(r)$, and state spaces $X(A)$ of substitutions or of subshifts.
%
Recall that the translations
by the integers $\mathbb{Z}$ on $L^2(\mathbb{R})$ are multiplication in the Fourier dual.
 For the Julia sets, this corresponds to coordinate
multiplication, i.e., multiplication by $z$ on the Hilbert space $L^2(X,\mu)$, consequently
a normal operator. The construction will involve Markov processes and
martingales.


Thus, many of the concepts related to multiresolutions in
$L^2(\mathbb{R})$ have analogues in a more general context; moreover, they
seem to exhibit interesting connections to the geometry of Julia sets.

\subsubsection{\label{Pyr}Pyramid algorithms and geometry}

Several of
the co-authors' recent projects
involve some
variant or other of the Ruelle transfer operator $R$, also called the
Ruelle--Perron--Frobenius operator. In each application, it arises in a
wavelet-like setting. But the settings are different from the familiar
$L^{2}(\mathbb{R})$-wavelet setup: one is for affine fractals, and the other
for Julia sets generated from iteration dynamics of rational functions $r(z)$.

Thus, there are two general themes in this work. In rough outline, they are
like this:
(1)~In the paper \cite{DuJo03}, Dutkay and Jorgensen construct a new class of
wavelets on extended fractals based on Cantor-like iterated function
systems, e.g., the middle-third Cantor set. Our Hilbert space in each of our
constructions is separable. For the fractals, it is built on Hausdorff measure
of the same (fractal) dimension as the IFS-Cantor set in question. In the
paper \cite{DuJo03}, we further introduce an associated Ruelle operator $R$,
and it plays a central role. For our wavelet construction, there is, as is to
be expected, a natural (equilibrium) measure $\nu$ which satisfies $\nu
R=\nu$, i.e., a left Perron--Frobenius eigenvector. It corresponds to the
Dirac point-measure on the frequency variable $\omega=0$ (i.e., low-pass) in
the standard $L^{2}(\mathbb{R})$-wavelet setting. It turns out that our
measures $\nu$ for the IFS-Cantor sets are given by infinite Riesz products,
and they are typically singular and have support${}={}$the whole circle
$\mathbb{T}$. This part of our research is related to recent, completely
independent, work by Benedetto et al.\ on Riesz products.

(2)~A second related
research direction is a joint project in progress with Ola Bratteli, where we
build wavelets on Hilbert spaces induced by Julia sets of rational functions
$r(z)$ of one complex variable. Let $r$ be a rational function of degree at
least $2$, and let $X(r)$ be the corresponding Julia set. Then there is a
family of Ruelle operators indexed by certain potentials, and a corresponding
family of measures $\nu$ on $X(r)$ which satisfy $\nu R=\nu$, again with the
measures $\nu$ playing a role somewhat analogous to the Dirac point-mass on
the frequency variable $\omega=0$, for the familiar $L^{2}(\mathbb{R})$-MRA
wavelet construction.

Independently of Dutkay--Jorgensen \cite{DuJo03},
John Benedetto and his co-authors Erica Bernstein
and Ioannis Konstantinidis \cite{BBK04}
have developed a new Fourier/in\-fi\-nite-product
approach to the very same singular measures
that arise in the study in \cite{DuJo03}.
The motivations and applications are different.
In \cite{DuJo03} the issue is wavelets on fractals,
and in \cite{BBK04}, it is time-frequency duality.

In this second class of problems, the measures $\nu$ are far less well
understood; and yet, if known, they would yield valuable insight into the
spectral theory of $r$-iteration systems and the corresponding Julia sets
$X(r)$.

\subsection{\label{Julia}Julia sets from complex dynamics}

Wavelet-like constructions are just under the surface in the following brief
sketch of problems. The projects all have to do with the kind of dynamical
systems $X$ already mentioned; and $\Omega$ is one of the standard probability
spaces. The space $X$ could be the circle, or a torus, or a solenoid, or a
fractal, or it could be a Julia set defined from a rational map in one complex
variable. More generally, consider maps $r$ of $X$ \emph{onto} $X$ which
generate some kind of dynamics. The discussion is restricted to the case when
for each $x$ in $X$, the pre-image $r^{-1}(\{x\})$ is assumed finite. Using
then a construction of Kolmogorov, it can be seen that these systems $X$ admit
useful measures which are obtained by a limit construction and by averaging
over the finite sets $r^{-1}(\{x\})$; see, e.g., Jorgensen's new monograph \cite{Jor04c} and Proposition \ref{prop3_1} for more details. We list three related problems.

(1)~{\bf Operator algebras}: A more systematic approach, along
these lines, to crossed products with endomorphisms. For example, trying to
capture the framework of en\-do\-mor\-phism-crossed products introduced
earlier by Bost and Connes \cite{BoCo95}. (2)~{\bf Dynamics}: Generalization to the study of onto-maps, when
the number and nature of branches of the inverse is more complicated; i.e.,
where it varies over $X$, and where there might be overlaps. (3)~{\bf Geometry}: Use of more detailed geometry and topology of Julia sets,
in the study of the representations that come from multiplicity
considerations. The Dutkay--Jorgensen paper \cite{DuJo03} on fractals is
actually a special case of Julia-set analysis. (One of the Julia sets $X$ from
the above \emph{is} a Cantor set in the usual sense, e.g., the Julia set of
$r(z) = 2z - 1/z$; or $r(z) = z^{2} - 3$.)

A main theme will be Hilbert spaces built on
Julia sets $X\left(r\right)$ built in turn on
rational functions $r\left(z\right)$ in
one complex variable. In the case when
$r$ is not a monomial, say of degree
$N \ge2$, we construct
$N$ functions $m_{i}$, $i=0,1,\dots,N-1$, on
$X\left(r\right)$ which satisfy a set of quadratic
conditions (\ref{eq2}) analogous to the axioms that define wavelet filters. We
have done this, so far, using representation theory, in special cases. But our
results are not yet in a form where they can be used in computations. By using results on the Hausdorff dimension
of Julia sets, it would also be interesting to identify inequivalent
representations, both in the case when the Julia set $X$ is fixed but the
functions vary, and in the case when $X$ varies.

Moreover, we  get an infinite-dimensional
\textquotedblleft loop group\textquotedblright\ $G=G(X,N)$ acting transitively
on these sets of functions $m_{i}$. The group $G$ consists of measurable
functions from $X$ into the group $\mathrm{U}_{N}$ of all unitary $N$ by $N$
complex matrices. In particular, our method yields information about this
group. Since $X$ can be geometrically complicated, it is not at all clear how
such matrix functions $X\rightarrow\mathrm{U}_{N}(\mathbb{C})$ might be
constructed directly.

The group $G$ consists of measurable functions $A\colon X(r)\rightarrow
\mathrm{U}_{N}(\mathbb{C})$, and the action of $A$ on an $(m_{i})$-system is
as follows: $(m_{i})\mapsto(m_{i}^{(A)})$, where
\begin{equation}
m_{i}^{(A)}(x)=\sum_{j}A_{i,j}(r(x))m_{j}(x),\qquad x\in X(r).
\label{miAxa}
\end{equation}

Even when $N=2$, the simplest non-monomial examples
include the Julia sets of $r(z)=z^{2}+c$, and the two functions $m_{0}$ and
$m_{1}$ on $X(c)=\operatorname{Julia}(z^{2}+c)$ are not readily available by
elementary means.
The use of operator algebras, representation theory seems natural for this
circle of problems.

\section{\label{ProPro}Main results}

\subsection{\label{MRonJ}Multiresolution/wavelet analysis on Julia sets}

We attempt to use equilibrium measures
(Brolin \cite{MR33:2805},
Ruelle \cite{Rue78a},
Mauldin--Urbanski \cite{MR2003772})
from discrete dynamical systems to specify generalized low-pass conditions in the Julia-set theory.  This would seem to be natural, but there are essential difficulties to overcome.
While Ruelle's thermodynamical formalism was developed for transfer operators where the weight function $W$ is strictly positive, the applications we have in mind dictate careful attention to the case when $W$ is not assumed strictly positive. The reason is that for our generalized multiresolution analysis on Julia sets, $W$ will be the absolute square of the low-pass filter.

\par
We begin with one example and a lemma. The example realizes
a certain class of state spaces from symbolic dynamics; and the lemma
spells out a class of invariant measures which will be needed in our
Hilbert-space constructions further into the section. The measures from
part (\ref{lemfme3_2(2)}) in Lemma \ref{lemfme3_2} below are often called strongly invariant
measures, or equilibrium measures from Ruelle's thermodynamical
formalism, see \cite{Rue78a}.

\begin{example}\label{exfme3}
Let $N\in\bz_+$, $N\geq2$ and let $A=(a_{ij})_{i,j=1}^N$ be an $N$
by $N$ matrix with all $a_{ij}\in\{0,1\}$. Set
$$X(A):=\{\,(x_i)\in\prod_{\mathbb{N}}\{1,\dots,N\}\mid A(x_i,x_{i+1})=1\,\}$$
and let $r=r_A$ be the restriction of the shift to $X(A)$, i.e.,
$$r_A(x_1,x_2,\dots)=(x_2,x_3,\dots),\qquad x=(x_1,x_2,\dots)\in X(A).$$
\begin{lemma}\label{lemfme3_1}
Let $A$ be as above. Then
$$\#\,r_A^{-1}(x)=\#\,\{\,y\in\{1,\dots,N\}\mid A(y,x_1)=1\,\}.$$
\end{lemma}
\par
It follows that $r_A \colon X(A)\rightarrow X(A)$ is onto iff $A$ is {\it
irreducible}, i.e., iff for all $j\in\bz_N$, there exists an
$i\in\bz_N$ such that $A(i,j)=1$. \par Suppose in addition that
$A$ is {\it aperiodic}, i.e., there exists $p\in\bz_+$ such that
$A^p>0$ on $\bz_N\times\bz_N$. We have the following lemma.
\begin{lemma}\label{lemfme3_2}\textup{(D. Ruelle
\cite{Rue89,Bal00})}\enspace Let $A$ be irreducible and aperiodic
and let $\phi\in C(X(A))$ be given. Assume that $\phi$ is a
Lipschitz function.
\begin{enumerate}
\item\label{lemfme3_2(1)}Set
$$(R_\phi f)(x)=\sum_{r_A(y)=x}e^{\phi(y)}f(y),\mbox{\quad for }f\in
C(X(A)).$$ Then there exist $\lambda_0>0$,
$$\lambda_0=\sup\{\,|\lambda|\mid\lambda\in\operatorname{spec}(R_\phi)\,\},$$
$h\in C(X(A))$ strictly positive and $\nu$ a Borel measure on
$X(A)$ such that
$$R_\phi h=\lambda_0 h,$$
$$\nu R_\phi=\lambda_0\nu,$$
and $\nu(h)=1$. The data are unique. \item\label{lemfme3_2(2)} In particular, setting
$$(R_0f)(x)=\frac{1}{\#\,r_A^{-1}(x)}\sum_{r_A(y)=x}f(y),$$
we may take $\lambda_0=1$, $h=1$ and $\nu=:\mu_A$, where $\mu_A$
is a probability measure on $X(A)$ satisfying the strong
invariance property
$$\int_{X(A)}f\,d\mu_A=\int_{X(A)}\frac{1}{\#\,r_A^{-1}(x)}\sum_{r_A(y)=x}f(y)\,d\mu_A(x),\qquad f\in
L^\infty(X(A).$$
\end{enumerate}
\end{lemma}
\end{example}
\par
Our next main results, Theorems \ref{th2_3} and Theorem \ref{th2_5}, are
operator-theoretic, and they address the following question: Starting
with a state space $X$ for a dynamical system, what are the conditions
needed for building a global Hilbert space $H(X)$ and an associated
multiscale structure? Our multiscales will be formulated in the context
of Definitions \ref{defhomog} (or equivalently \ref{defmulti}) above; that is, we define our
multiscale structure in terms of a natural Hilbert space $H(X)$
and a
certain unitary operator $U$ on $H(X)$ which implements the multiscales.
\par
 After our presentation of the Hilbert-space context, we turn to
spectral theory: In our next result, we revisit Baggett's
dimension-consistency equation. The theorem of Baggett et al.\ \cite{BJMP04} concerns
a condition on a certain pair of multiplicity functions for the existence
of a class of wavelet bases for $L^2(\Br^d)$. In Theorem \ref{thmul1} below, we
formulate the corresponding multiplicity functions in the context of
nonlinear dynamics, and we generalize the Baggett et al.\ theorem to the
Hilbert spaces $H(X)$ of our corresponding state spaces $X$.
\par As preparation for Theorem \ref{thpr_5} below, we introduce a family of
$L^2$-mar\-tin\-gales, and we prove that the Hilbert spaces $H(X)$ are $L^2$-mar\-tin\-gale Hilbert spaces. In this context, we prove the existence
of an explicit wavelet transform (Theorem \ref{thpr_5}) for our nonlinear state
spaces. But in the martingale context, we must work with dilated state
spaces $X_\infty$; and in Proposition \ref{prop3_1}, and Theorems \ref{th4_18}, \ref{th4_19}, \ref{th4_14}, and \ref{th4_17}, we outline the role of our multiscale structure in
connection with certain random-walk measures on $X_\infty$. 
A second issue concerns extreme measures: in Theorem \ref{th7_2} 
we characterize the extreme points in the convex set of these
random-walk measures.
\begin{theorem}\label{th2_3}
Let $\mathcal{A}$ be a unital $C^*$-algebra, $\alpha$ an
endomorphism on $\mathcal{A}$, $\mu$ a state on $\mathcal{A}$ and,
$m_0\in\mathcal{A}$, such that
\begin{equation}\label{eq2_5}
\mu(m_0^*\alpha(f)m_0)=\mu(f),\qquad f\in\mathcal{A}.
\end{equation}
Then there exists a Hilbert space $H$, a representation $\pi$ of
$\mathcal{A}$ on $H$, $U$ a unitary on $H$, and a vector
$\varphi\in\mathcal{A}$, with the following properties:
\begin{equation}\label{eq2_6}
U\pi(f)U^*=\pi(\alpha(f)),\qquad f\in\mathcal{A},
\end{equation}
\begin{equation}\label{eq2_7}
\ip{\varphi}{\pi(f)\varphi}=\mu(f),\qquad f\in\mathcal{A},
\end{equation}
\begin{equation}\label{eq2_8}
U\varphi=\pi(\alpha(1)m_0)\varphi
\end{equation}
\begin{equation}\label{eq2_9}
\cj{\operatorname{span}}\{\,U^{-n}\pi(f)\varphi\mid n\geq0,f\in\mathcal{A}\,\}=H.
\end{equation}
Moreover, this is unique up to an intertwining isomorphism.
\end{theorem}
We call $(H,U,\pi,\varphi)$ the covariant system associated to
$\mu$ and $m_0$.

\begin{corollary}\label{cor2_4}
Let $X$ be a measure space, $r \colon X\rightarrow X$ a measurable, onto
map and $\mu$ a probability measure on $X$ such that
\begin{equation}
\int_Xf\,d\mu=\int_X\frac{1}{\#\,r^{-1}(x)}\sum_{r(y)=x}f(y)\,d\mu(x).
\end{equation}
Let $h\in L^1(X)$, $h\geq0$ such that $$\frac{1}{\#\,
r^{-1}(x)}\sum_{r(y)=x}|m_0(y)|^2h(y)=h(x),\qquad x\in X.$$ Then
there exists \textup{(}uniquely up to isomorphism\/\textup{)} a Hilbert space $H$, a
unitary $U$, a representation $\pi$ of $L^\infty(X)$ and a vector
$\varphi\in H$ such that
$$U\pi(f)U^{-1}=\pi(f\circ r),\qquad f\in L^\infty(X),$$
$$\ip{\varphi}{\pi(f)\varphi}=\int_Xfh\,d\mu,\qquad f\in
L^\infty(X),$$
$$U\varphi=\pi(m_0)\varphi,$$
$$\cj{\operatorname{span}}\{\,U^{-n}\pi(f)\varphi\mid n\geq0,f\in
L^\infty(X)\,\}=H.$$ We call $(H,U,\pi,\varphi)$ the covariant
system associated to $m_0$ and $h$.
\end{corollary}

\subsection{\label{BS}The Baumslag--Solitar group}

Our next theorem is motivated by the theory of representations
of a certain group which is called the Baumslag--Solitar group, and which
arises by the following simple semidirect product construction:
 the Baumslag--Solitar group $BS(1,N)$ is the group with two
generators $u$ and $t$ and one relation $utu^{-1}=t^N$, where $N$
is a fixed positive integer.
 Therefore the two operators $U$ of $N$-adic dilation and $T$ of integral translation operators on $\mathbb{R}$ give a representation of this group. But so do all familiar MRA-wavelet constructions.
\par
Representations of the Baumslag--Solitar group that admit wavelet
bases can be constructed also on some other spaces, such as
$\ltwor\oplus\dots\oplus\ltwor$ \cite{BDP04,Dut04a} or
some fractal spaces \cite{DuJo03}.
\par
     Hence wavelet theory fits naturally into this abstract
     setting of special representations of the group $BS(1,N)$ realized on
     some Hilbert space.

\par
It is already known \cite{Dut04d} that the representations that
admit wavelets are faithful and weakly equivalent to the right
regular representation of the group.
\par

\begin{theorem}\label{th2_5}
\textup{(i)} Let $H$ be a Hilbert space, $S$ an isometry on $H$. Then there
exist a Hilbert space $\hat H$ containing $H$ and a unitary $\hat
S$ on $\hat H$ such that
\begin{equation}\label{eq2_5_1}
\hat S|_{H}=S,
\end{equation}
\begin{equation}\label{eq2_5_2}
\cj{\bigcup_{n\geq0}\hat S^{-n}H}=\hat H.
\end{equation}
Moreover these are unique up to an intertwining isomorphism. \par
\textup{(ii)} If $\mathcal{A}$ is a $C^*$-algebra, $\alpha$ is an
endomorphism on $\mathcal{A}$ and $\pi$ is a representation of
$\mathcal{A}$ on $H$ such that
\begin{equation}\label{eq2_5_3}
S\pi(g)=\pi(\alpha(g))S,\qquad g\in\mathcal{A},
\end{equation}
then there exists a unique representation $\hat\pi$ on $\hat H$
such that
\begin{equation}\label{eq2_5_4}
\hat\pi(g)|_H=\pi(g),\qquad g\in\mathcal{A},
\end{equation}
\begin{equation}\label{eq2_5_5}
\hat S\hat\pi(g)=\hat\pi(\alpha(g))\hat S,\qquad g\in\mathcal{A}.
\end{equation}
\end{theorem}

\begin{corollary}\label{cor2_6}
Let $X,r,$ and $\mu$ be as in Corollary \ref{cor2_4}. Let $I$ be a
finite or countable set. Suppose
$H \colon X\rightarrow\mathcal{B}(\ell^2(I))$ has the property that
$H(x)\geq 0$ for almost every $x\in X$, and $H_{ij}\in L^1(X)$ for
all $i,j\in I$. Let $M_0 \colon X\rightarrow\mathcal{B}(\ell^2(I))$ such
that $x\mapsto\|M_0(x)\|$ is essentially bounded. Assume in
addition that
\begin{equation}\label{eq2_6_1}
\frac{1}{\#\,r^{-1}(x)}\sum_{r(y)=x}M_0^*(y)H(y)M_0(y)=H(x),\mbox{\quad for a.e. }x\in X.
\end{equation}
Then there exists a Hilbert space $\hat K$, a unitary operator
$\hat U$ on $\hat K$, a representation $\hat\pi$ of $L^\infty(X)$
on $\hat K$, and a family of vectors $(\varphi_i)\in\hat K$, such
that
$$\hat U\hat\pi(g)\hat U^{-1}=\hat\pi(g\circ r),\qquad g\in L^\infty(X),$$
$$\hat U\varphi_i=\sum_{j\in I}\hat\pi((M_0)_{ji})\varphi_j,\qquad i\in I,$$
$$\ip{\varphi_i}{\hat\pi(f)\varphi_j}=\int_XfH_{ij}\,d\mu,\qquad i,j\in I,\;f\in L^\infty(X),$$
$$\cj{\operatorname{span}}\{\,\hat\pi(f)\varphi_i\mid n\geq0,f\in L^\infty(X),i\in I\,\}=\hat K.$$
These are unique up to an intertwining unitary isomorphism.
\textup{(}All functions are assumed weakly measurable in the sense that $x\mapsto\ip{\xi}{F(x)\eta}$ is measurable for all $\xi,\eta\in \ell^2(I)$.\textup{)}
\end{corollary}

\par
Suppose now that $H$ is a Hilbert space with an isometry $S$ on it
and with a normal representation $\pi$ of $L^\infty(X)$ on $H$
that satisfies the covariance relation
\begin{equation}\label{eqmul2}
S\pi(g)=\pi(g\circ r)S,\qquad g\in L^\infty(X).
\end{equation}
\par
Theorem \ref{th2_5} shows that there exists a Hilbert space $\hat
H$ containing $H$, a unitary $\hat S$ on $\hat H$  and a
representation $\hat\pi$ of $L^\infty(X)$ on $\hat H$ such that
$$(V_n:=\hat S^{-n}(H))_n\mbox{ form an increasing sequence of
subspaces with dense union},$$
$$\hat S|_H=S,$$
$$\hat\pi|_H=\pi,$$
$$\hat S\hat\pi(g)=\hat\pi(g\circ r)\hat S.$$

\begin{theorem}\label{thmul1}
\textup{(i)} $V_1=\hat S^{-1}(H)$ is invariant for the representation
$\hat\pi$. The multiplicity functions of the representation
$\hat\pi$ on $V_1$, and on $V_0=H$, are related by
\begin{equation}\label{eqmul3}
m_{V_1}(x)=\sum_{r(y)=x}m_{V_0}(y),\qquad x\in X.
\end{equation}
\par
\textup{(ii)} If $W_0:=V_1\ominus V_0=\hat S^{-1}H\ominus H$, then
\begin{equation}\label{eqmul4}
m_{V_0}(x)+m_{W_0}(x)=\sum_{r(y)=x}m_{V_0}(y),\qquad x\in X.
\end{equation}
\end{theorem}
\begin{proof}
Note that $\hat S$ maps $V_1$ to $V_0$, and the covariance
relation implies that the representation $\hat\pi$ on $V_1$ is
isomorphic to the representation $\pi^r \colon g\mapsto\pi(g\circ r)$ on
$V_0$. Therefore we have to compute the multiplicity of the
latter, which we denote by $m^r_{V_0}$.
\par
By the spectral theorem there exists a unitary isomorphism
$J \colon H(=V_0)\rightarrow L^2(X,m_{V_0},\mu)$, where, for a {\it
multiplicity function} $m \colon X\rightarrow\{0,1,\dots,\infty\}$, we use
the notation
$$
L^2(X,m,\mu):=\{\,f \colon X\rightarrow\cup_{x\in X}\bc^{m(x)}\mid
f(x)\in\bc^{m(x)},\int_X\|f(x)\|^2\,d\mu(x)<\infty\,\}.$$ In
addition $J$ intertwines $\pi$ with the representation of
$L^\infty(X)$ by multiplication operators, i.e.,
$$(J\pi(g)J^{-1}(f))(x)=g(x)f(x),\qquad g\in L^\infty(X),f\in
L^2(X,m_{V_0},\mu),x\in X.$$

\begin{remark}\label{rem4_9}
Here we are identifying $H$ with $L^2(X,m_{V_0},\mu)$ via the {\it
spectral representation}. We recall the details of this
representation $H\ni f\mapsto\tilde f\in L^2(X,m_{V_0},\mu)$.
\par
Recall that any normal representation $\pi\in \operatorname{Rep}(L^\infty(X),H)$
is the orthogonal sum
\begin{equation}\label{eqmulrem*}
H=\sideset{}{^{\smash{\oplus}}}{\sum}\limits_{k\in C}[\pi(L^\infty(X))k],
\end{equation}
where the set $C$ of vectors $k\in H$ is chosen such that
\begin{itemize}
\item \makebox[\displayboxwidth]{$\displaystyle\|k\|=1,$}
\item \makebox[\displayboxwidth]{$\displaystyle
\ip{k}{\pi(g)k}=\int_Xg(x)v_k(x)^2\,d\mu(x),\mbox{ for all }k\in
C;
$}
\item \makebox[\displayboxwidth]{$\displaystyle
\ip{k'}{\pi(g)k}=0,\quad g\in L^\infty(X),k,k'\in C,k\neq
k';\mbox{ orthogonality}.$}
\end{itemize}
\par
The formula (\ref{eqmulrem*}) is obtained by a use of Zorn's
lemma. Here, $v_k^2$ is the Radon-Nikodym derivative of
$\ip{k}{\pi(\cdot)k}$ with respect to $\mu$, and we use that $\pi$
is assumed normal.
\par
For $f\in H$, set
$$f=\sideset{}{^{\smash{\oplus}}}{\sum}\limits_{k\in C}\pi(g_k)k,\quad g_k\in L^\infty(X)$$
and
$$\tilde f=\sideset{}{^{\smash{\oplus}}}{\sum}\limits_{k\in C}g_kv_k\in L_\mu^2(X,\ell^2(C)).$$
Then $Wf=\tilde f$ is the desired spectral transform, i.e.,
$$W\mbox{ is unitary},$$
$$W\pi(g)=M(g)W,$$
and
$$\|\tilde f(x)\|^2=\sum_{k\in C}|g_k(x)v_k(x)|^2.$$
Indeed, we have
\begin{multline*}
\int_X\|\tilde f(x)\|^2\,d\mu(x)=\int_X\sum_{k\in
C}|g_k(x)|^2v_k(x)^2\,d\mu(x)=\sum_{k\in
C}\int_X|g_k|^2v_k^2\,d\mu
\\
{}=\sum_{k\in C}\ip{k}{\pi(|g_k|^2)k}=\sum_{k\in
C}\|\pi(g_k)k\|^2=\left\|\sideset{}{^{\smash{\oplus}}}{\sum}\limits_{k\in C}
\pi(g_k)k\right\|^2_H=\|f\|^2_{H}.
\end{multline*}
It follows in particular that the multiplicity function
$m(x)=m_{H}(x)$ is
$$m(x)=\#\,\{\,k\in C\mid v_k(x)\neq0\,\}.$$
Setting
$$X_i:=\{\,x\in X\mid m(x)\geq i\,\},\qquad i\geq 1,$$
we see that
$$H\simeq\sideset{}{^{\smash{\oplus}}}{\sum}L^2(X_i,\mu)\simeq L^2(X,m,\mu),$$
and the isomorphism intertwines $\pi(g)$ with multiplication
operators.
\end{remark}
\par
Returning to the proof of the theorem, we have to find the similar
form for the representation $\pi^r$. Let
\begin{equation}\label{eqmul5}
\tilde m(x):=\sum_{r(y)=x}m_{V_0}(y),\qquad x\in X. \end{equation}
Define the following unitary isomorphism:
$$L \colon L^2(X,m_{V_0},\mu)\rightarrow L^2(X,\tilde m,\mu),$$
$$(L\xi)(x)=\frac{1}{\sqrt{\#\,r^{-1}(x)}}(\xi(y))_{r(y)=x}.$$
(Note that the dimensions of the vectors match because of
(\ref{eqmul5}).) This operator $L$ is unitary. For $\xi\in
L^2(X,m_{V_0},\mu)$, we have
\begin{align*}
\|L\xi\|^2_{L^2(X,m_{V_0},\mu)}&=\int_X\|L\xi(x)\|^2\,d\mu(x)\\
&=\int_X\frac{1}{\#\,r^{-1}(x)}\sum_{r(y)=x}\|\xi(y)\|^2\,d\mu(x)\\
&=\int_X\|\xi(x)\|^2\,d\mu(x).
\end{align*}

And $L$ intertwines the representations. Indeed, for $g\in
L^\infty(X)$,
$$L(g\circ r\,\xi)(x)=(g(r(y))\xi(y))_{r(y)=x}=g(x)L(\xi)(x).$$
Therefore, the multiplicity of the representation
$\pi^r \colon g\mapsto\pi(g\circ r)$ on $V_0$ is $\tilde m$, and this
proves (i).
\par
(ii) follows from (i).
\par
{\it Conclusions.} By definition, if $k\in C$,
$$\ip{k}{\pi(g)k}=\int_Xg(x)v_k(x)^2\,d\mu(x),\mbox{ and }$$
$$\ip{k}{\pi^r(g)k}=\int_Xg(r(x))v_k(x)^2\,d\mu(x)=\int_Xg(x)\frac{1}{\#\,r^{-1}(x)}\sum_{r(y)=x}v_k(x)^2\,d\mu(x);$$
and so
\begin{align*}
m^r(x)&=\#\,\{\,k\in
C\mid\sum_{r(y)=x}v_k(y)^2>0\,\}\\
&=\sum_{r(y)=x}\#\,\{\,k\in C\mid v_k(y)^2>0\,\}\\
&=\sum_{r(y)=x}m(y).
\end{align*}
Let $C^m(x):=\{\,k\in C\mid v_k(x)\neq0\,\}$. Then we showed that
$$C^m(x)=\bigcup_{y\in X, r(y)=x}C^m(y)$$
and that $C^m(y)\cap C^m(y')=\emptyset$ when $y\neq y'$ and
$r(y)=r(y')=x.$ Setting $\mathcal{H}(x)=\ell^2(C^m(x))$, we have
$$\mathcal{H}(x)=\ell^2(C^m(x))=\sideset{}{^{\smash{\oplus}}}{\sum}\limits_{r(y)=x}
\ell^2(C^m(y))=\sideset{}{^{\smash{\oplus}}}{\sum}\limits_{r(y)=x}\mathcal{H}(y).\eqno\qedhere$$
\end{proof}

\subsection{\label{projlim}Spectral decomposition of covariant representations: projective limits}
We give now a different representation of the construction of the
covariant system associated to $m_0$ and $h$.

We work in either the category of measure spaces or topological
spaces. \begin{definition}\label{def4_10} Let $r \colon X\rightarrow X$
be onto, and assume that $\#\, r^{-1}(x)<\infty$ for all $x\in X$.
We define the {\it projective limit} of the system:
\begin{equation}\label{eqpr_1}
X\stackrel{r}{\longleftarrow}
X\stackrel{r}{\longleftarrow}X\cdots\longleftarrow X_\infty
\end{equation} as
$$X_\infty:=\{\,\hat x=(x_0,x_1,\dots)\mid r(x_{n+1})=x_n,\mbox{ for all }n\geq0\,\}.$$
\end{definition}
\par Let $\theta_n \colon X_\infty\rightarrow X$ be the projection
onto the $n$-th component:
$$\theta_n(x_0,x_1,\dots)=x_n,\qquad(x_0,x_1,\dots)\in X_\infty.$$
Taking inverse images of sets in $X$ through these projections, we
obtain a sigma-algebra on $X_\infty$, or a topology on $X_\infty$.
\par We have an induced mapping $\hat r \colon  X_\infty\rightarrow
X_\infty$ defined by
\begin{equation}\label{eqpr_3}
\hat r(\hat x)=(r(x_0),x_0,x_1,\dots),\mbox{ and with inverse }\hat
r^{-1}(\hat x)=(x_1,x_2,\dots). \end{equation} so $\hat r$ is an
automorphism, i.e., $\hat r\circ\hat r^{-1}=\mbox{id}_{X_\infty}$
and $\hat r^{-1}\circ\hat r=\mbox{id}_{X_\infty}$.
\par
Note that
$$\theta_n\circ\hat r=r\circ \theta_n=\theta_{n-1},$$
$$
\setlength{\unitlength}{0.075\textwidth}
\begin{picture}(4.8,3.5)(-0.9,-0.75)
\put(0.075,2.25){\vector(1,0){2.85}}
\put(0,2){\vector(3,-2){3}}
\put(3.375,1.925){\vector(0,-1){1.85}}
\put(-0.375,2.25){\makebox(0,0){$X_\infty$}}
\put(3.375,2.25){\makebox(0,0){$X$}}
\put(3.375,-0.25){\makebox(0,0){$X$}}
\put(1.5,2.375){\makebox(0,0)[b]{$\theta_n$}}
\put(1.5,1){\makebox(0,0)[tr]{$\theta_{n-1}$}}
\put(3.5,1){\makebox(0,0)[l]{$r$}}
\end{picture}
\begin{picture}(4.8,3.5)(-0.9,-0.75)
\put(0.075,2.25){\vector(1,0){2.85}}
\put(0,2){\vector(3,-2){3}}
\put(3.375,1.925){\vector(0,-1){1.85}}
\put(-0.375,2.25){\makebox(0,0){$X_\infty$}}
\put(3.375,2.25){\makebox(0,0){$X_\infty$}}
\put(3.375,-0.25){\makebox(0,0){$X$}}
\put(1.5,2.375){\makebox(0,0)[b]{$\hat r$}}
\put(1.5,1){\makebox(0,0)[tr]{$\theta_{n-1}$}}
\put(3.5,1){\makebox(0,0)[l]{$\theta_n$}}
\end{picture}
\begin{picture}(3.7,3.5)(-0.85,-0.75)
\put(0.075,2.25){\vector(1,0){1.85}}
\put(2.375,1.925){\vector(0,-1){1.85}}
\put(-0.375,1.925){\vector(0,-1){1.85}}
\put(0.075,-0.25){\vector(1,0){1.85}}
\put(-0.375,2.25){\makebox(0,0){$X_\infty$}}
\put(2.375,2.25){\makebox(0,0){$X_\infty$}}
\put(-0.375,-0.25){\makebox(0,0){$X$}}
\put(2.375,-0.25){\makebox(0,0){$X$.}}
\put(1,2.375){\makebox(0,0)[b]{$\hat r$}}
\put(-0.25,1){\makebox(0,0)[l]{$\theta_n$}}
\put(1,-0.375){\makebox(0,0)[t]{$r$}}
\put(2.5,1){\makebox(0,0)[l]{$\theta_n$}}
\end{picture}
$$
\par
Consider a probability measure $\mu$ on $X$ that satisfies
\begin{equation}\label{eqpr_4}
\int_Xf\,d\mu=\int_X\frac{1}{\#\,r^{-1}(x)}\sum_{r(y)=x}f(y)\,d\mu(x).
\end{equation}
\par

For $m_0\in L^\infty(X)$, define
\begin{equation}\label{eqpr_7}
(R\xi)(x)=\frac{1}{\#\,
r^{-1}(x)}\sum_{r(y)=x}|m_0(y)|^2\xi(y),\qquad \xi\in L^1(X).
\end{equation}
\par
         We now resume the operator-theoretic theme of our paper, that of
extending a system of operators in a given Hilbert space to a ``dilated''
system in an ambient or extended Hilbert space, the general idea being that
the dilated system acquires a ``better'' spectral theory: for example,
contractive operators dilate to unitary operators in the extended Hilbert
space, and similarly, endomorphisms dilate to automorphisms. This, of course,
is a general theme in operator theory; see, e.g., \cite{BaVi02} and \cite{MR2004i:46109}. But in our
present setting, there is much more structure than the mere Hilbert-space
geometry. We must adapt the underlying operator theory to the particular
function spaces and measures at hand. The next theorem (Theorem \ref{thpr_1}) is key
to the understanding of the dilation step as it arises in our context of
multiscale theory. It takes into consideration how we arrive at our function
spaces by making dynamics out of the multiscale setting. Hence, the
dilations we construct take place at three separate levels, as follows:
\begin{itemize}
\item dynamical systems,
$$(X,r,\mu)\mbox{ endomorphism }\rightarrow(X_\infty,\hat
r,\hat\mu),\mbox{ automorphism };$$ \item Hilbert spaces,
$$L_2(X,h\,d\mu)\rightarrow (R_{m_0}h=h)\rightarrow
L^2(X_\infty,\hat\mu);$$ \item operators,
$$S_{m_0}\mbox{ isometry }\rightarrow U\mbox{ unitary (if
}m_0\mbox{ is non-sin\-gu\-lar}),$$
$$M(g)\mbox{ multiplication operator }\rightarrow M_\infty(g).$$
\end{itemize}
\begin{definition}
A function $m_0$ on a measure space is called {\it singular} if
$m_0=0$ on a set of positive measure.
\end{definition}
\par
In general, the operators $S_{m_0}$ on $H_0=L^2(X,h\,d\mu)$, and
$U$ on $L^2(X_\infty,\hat\mu)$, may be given only by abstract
Hilbert-space axioms; but in our {\it martingale representation},
we get the following two concrete formulas:
$$(S_{m_0}\xi)(x)=m_0(x)\xi(r(x)),\qquad x\in X,\xi\in H_0,$$
$$(Uf)(\hat x)=m_0(x_0)f(\hat r(\hat x)),\qquad\hat x\in
X_\infty,f\in L^2(X_\infty,\hat\mu).$$
\begin{definition} Let $X$, $r$, $\mu$, and $W$ be given as before. Recall $r  \colon 
X \rightarrow X$ is a finite-to-one map onto $X$, $\mu$ is a
corresponding strongly invariant measure, and $W  \colon  X \rightarrow
[0, \infty)$ is a measurable weight function. We define the
corresponding Ruelle operator $R = R_W$, and we let $h$ be a
chosen Perron--Frobenius--Ruelle eigenfunction. From these data, we
define the following sequence of measures $\omega_n$ on $X$. (Our
presentation of a wavelet transform in Theorem \ref{thpr_5} will
depend on this sequence.) For each $n$, the measure $\omega_n$ is
defined by the following formula:
\begin{equation}\label{eqpr_6}
\omega_n(f)=\int_XR^n(fh)\,d\mu,\qquad f\in L^\infty(X).
\end{equation}
\end{definition}

\par
To build our wavelet transform, we first prove in Theorem
\ref{thpr_1} below that each of the measure families $(\omega_n)$
defines a unique $W$-quasi-invariant measure $\hat\mu$ on
$X_\infty$. In the theorem, we take $W := |m_0 |^2$, and we show
that $\hat\mu$ is quasi-invariant with respect to $\hat r$ with
transformation Radon-Nikodym derivative $W$.

\begin{theorem}\label{thpr_1}
If $h\in L^1(X)$, $h\geq0$ and $Rh=h$, then there exists a unique
measure $\hat\mu$ on $X_\infty$ such that
$$\hat\mu\circ\theta_n^{-1}=\omega_n,\qquad n\geq0.$$

We can identify functions on $X$ with functions on
$X_\infty$ by
$$f(x_0,x_1,\dots)=f(x_0),\qquad f \colon X\rightarrow\bc.$$
Under this identification,
\begin{equation}\label{eqpr_10}
\frac{d(\hat\mu\circ\hat r)}{d\hat\mu}=|m_0|^2.
\end{equation}
\end{theorem}

\begin{theorem}\label{thpr_3}
Suppose $m_0$ is non-sin\-gu\-lar, i.e., it does not vanish on a set
of positive measure. Define $U$ on $L^2(X_\infty,\hat\mu)$ by
$$Uf=m_0f\circ\hat r,\qquad f\in L^2(X_\infty,\hat\mu),$$
$$\pi(g)f=gf,\qquad g\in L^\infty(X),f\in L^2(X_\infty,\hat\mu),$$
$$\varphi=1.$$
Then $(L^2(X_\infty,\hat\mu),U,\pi,\varphi)$ is the covariant
system associated to $m_0$ and $h$ as in Corollary \ref{cor2_4}.
Moreover, if $M_gf=gf$ for $g\in L^\infty(X_\infty,\hat\mu)$ and
$f\in L^2(X_\infty,\hat\mu)$, then
$$UM_gU^{-1}=M_{g\circ\hat r}.$$
\end{theorem}

\par
The Hilbert space $L^2(X_\infty,\hat\mu)$ admits a different representation as an $L^2$-mar\-tin\-gale Hilbert space.
Let
$$H_n:=\{\,f\in L^2(X_\infty,\hat\mu)\mid
f=\xi\circ\theta_n,\xi\in L^2(X,\omega_n)\,\}.$$ Then $H_n$ form an
increasing sequence of closed subspaces which have dense union.
\par
We can identify the functions in $H_n$ with functions in
$L^2(X,\omega_n)$, by
$$i_n(\xi)=\xi\circ\theta_n,\qquad\xi\in L^2(X,\omega_n).$$
The definition of $\hat\mu$ makes $i_n$ an isomorphism between
$H_n$ and $L^2(X,\omega_n).$
\par
Define
$$\mathcal{H}:=\left\{(\xi_0,\xi_1,\dots)\Bigm| \xi_n\in
L^2(X,\omega_n),R(\xi_{n+1}h)=\xi_nh,\,
\sup_n\int_X\!R^n(|\xi_n|^2h)\,d\mu<\infty\right\}$$ with the scalar
product
$$\ip{(\xi_0,\xi_1,\dots)}{(\eta_0,\eta_1,\dots)}=\lim_{n\rightarrow\infty}\int_XR^n(\cj\xi_n\eta_nh)\,d\mu.$$
\begin{theorem}\label{thpr_5}
The map $\Phi \colon L^2(X_\infty,\hat\mu)\rightarrow\mathcal{H}$ defined
by
$$\Phi(f)=(i_n^{-1}(P_nf))_{n\geq0},$$
where $P_n$ is the projection onto $H_n$, is an isomorphism. The
transform $\Phi$ satisfies the following three conditions, and it
is determined uniquely by them:
$$\Phi U\Phi^{-1}(\xi_n)_{n\geq0}=(m_0\circ
r^n\,\xi_{n+1})_{n\geq0},$$
$$\Phi\pi(g)\Phi^{-1}(\xi_n)_{n\geq0}=(g\circ
r^n\,\xi_n)_{n\geq0},$$
$$\Phi\varphi=(1,1,\dots).$$
\end{theorem}

\begin{theorem}\label{thcond3}
There exists a unique isometry
$\Psi \colon L^2(X_\infty,\hat\mu)\rightarrow\tilde H$ such that
$$\Psi(\xi\circ\theta_n)=\tilde U^{-n}\tilde\pi(\xi)\tilde
U^n\tilde\varphi,\qquad\xi\in L^\infty(X,\mu).$$ $\Psi$
intertwines the two systems, i.e.,
$$\Psi U=\tilde U\Psi,\; \Psi\pi(g)=\tilde\pi(g)\Psi,\mbox{\quad for }g\in
L^\infty(X,\mu),\;\Psi\varphi=\tilde\varphi.$$
\end{theorem}

\begin{theorem}\label{corinter2}
Let $(m_0,h)$ be a Perron--Ruelle--Frobenius pair with $m_0$
non-sin\-gu\-lar.
\begin{enumerate}
\item\label{corinter2(1)} For each operator $A$ on $L^2(X_\infty,\hat\mu)$ which
commutes with $U$ and $\pi$, there exists a cocycle $f$, i.e., a
bounded measurable function $f \colon X_\infty\rightarrow\bc$ with
$f=f\circ\hat r$, $\hat\mu$-a.e., such that
\begin{equation}\label{eqinterc1}A=M_f,
\end{equation}
and, conversely each cocycle defines an operator in the commutant.
\item\label{corinter2(2)} For each measurable harmonic function $h_0 \colon X\rightarrow\bc$,
i.e., $R_{m_0}h_0=h_0$, with $|h_0|^2\leq ch^2$ for some $c\geq0$,
there exists a unique cocycle $f$ such that
\begin{equation}\label{eqinterc2}
h_0=E_0(f)h,
\end{equation}
where $E_0$ denotes the conditional expectation from $L^\infty(X_\infty,\hat\mu)$ to $L^\infty(X,\mu)$, and conversely, for each cocycle the function $h_0$ defined by
\textup{(\ref{eqinterc2})} is harmonic. \item\label{corinter2(3)} The correspondence
$h_0\rightarrow f$ in \textup{(\ref{corinter2(2)})} is given by
\begin{equation}\label{eqinterc3}
f=\lim_{n\rightarrow\infty}\frac{h_0}{h}\circ\theta_n
\end{equation}
where the limit is pointwise $\hat\mu$-a.e., and in
$L^p(X_\infty,\hat\mu)$ for all $1\leq p<\infty$.
\end{enumerate}
\end{theorem}
\par
The next result in this section concerns certain path-space measures,
indexed by a base space $X$. Such a family is also called a process: it is a
family of positive Radon measures $P_x$, indexed by $x$ in the base space $X$.
Each $P_x$ is a measure on the probability space $\Omega$ which is constructed
from $X$ by the usual Cartesian product, i.e., countably infinite Cartesian
product of $X$ with itself. The Borel structure on $\Omega$ is generated by the
usual cylinder subsets of $\Omega$. Given a weight function $W$ on $X$, the
existence of the measures $P_x$ comes from an application of a general
principle of Kolmogorov.

\par
Let $X$ be a metric space and $r \colon X\rightarrow X$ an $N$ to $1$
map. Denote by $\tau_k \colon X\rightarrow X$, $k\in\{1,\dots,N\}$, the
branches of $r$, i.e., $r(\tau_k(x))=x$ for $x\in X$, the sets
$\tau_k(X)$ are disjoint and they cover $X$.
\par
Let $\mu$ be a measure on $X$ with the property
\begin{equation}\label{eq3_1}
\mu=\frac{1}{N}\sum_{k=1}^N\mu\circ\tau_k^{-1}.
\end{equation}
This can be rewritten as
\begin{equation}\label{eq3_2}
\int_{X}f(x)\,d\mu(x)=\frac{1}{N}\sum_{k=1}^N\int_Xf(\tau_k(x))\,d\mu(x),
\end{equation}
which is equivalent also to the strong invariance property. \par

Let $W$, $h\geq0$ be two functions on $X$ such that
\begin{equation}\label{eq3_3}
\sum_{k=1}^NW(\tau_k(x))h(\tau_k(x))=h(x),\qquad x\in X.
\end{equation}
Denote by $\Omega$ the multi-index set
$$\Omega:=\Omega_N:=\prod_{\mathbb{N}}\{1,\dots,N\}.$$ Also we define
$$W^{(n)}(x):=W(x)W(r(x))\cdots W(r^{n-1}(x)),\qquad x\in X.$$

\par
We will return to a more general variant of these path-%
space measures in Theorem \ref{th4_17} in the next section. This
will be in the context of a system $(X, r, V)$ where $r  \colon  X
\rightarrow X$ is a give finite to one endomorphism, and $W  \colon  X
\rightarrow [0, \infty)$ is a measurable weight function. In this
context the measures $(P_x)$ are indexed by a base space $X$, and
they are measures on the projective space $X_\infty$ built on $(X,
r)$, see Definition \ref{def4_10}. The family $(P_x)$ is also
called a process: each $P_x$ is a positive Radon measure on
$X_\infty$, and $X_\infty$ may be thought of as a discrete
path space.
\par
The Borel structure on $\Omega$, and on $X_\infty$, is generated by
the usual cylinder subsets. Given a weight function $W$ on $X$,
the existence of the measures $P_x := P^W_x$ comes from an
application of a general principle of Kolmogorov.

\begin{proposition}\label{prop3_1}
For every $x\in X$ there exists a positive Radon measure $P_x$ on
$\Omega$ such that, if $f$ is a bounded measurable function on $X$
which depends only on the first $n$ coordinates
$\omega_1,\dots,\omega_n$, then
\begin{multline}\label{eq3_4}\int_\Omega
f(\omega)\,dP_x(\omega)\\=\!\sum_{\omega_1,\dots,\omega_n}\!W^{(n)}(\tau_{\omega_n}\tau_{\omega_{n-1}}\cdots\tau_{\omega_1}(x))h(\tau_{\omega_n}\tau_{\omega_{n-1}}\cdots\tau_{\omega_1}(x))f(\omega_1,\dots,\omega_n).
\end{multline}
\end{proposition}

\section{\label{other}Remarks on other applications}

\subsection{\label{Waveletsets}Wavelet sets}
Investigate the existence and construction of wavelet sets and elementary wavelets and frames in the Julia-set theory, and the corresponding interpolation theory, and their relationships to
generalized multiresolution analysis
(GMRA).  The unitary system approach to wavelets by Dai and Larson \cite{DL98} dealt with systems that can be very irregular. And recent work by \'Olafsson et al.\ shows that wavelet sets can be essential to a basic dilation-translation wavelet theory even for a system where the set of dilation unitaries is not a group.

\par

\subsection{\label{renormalization}The renormalization question}
When are there renormalizable iterates of $r$ in arithmetic progression, i.e., when are there iteration periods $n$ such that the system $\{r^{kn}, k \in \mathbb{N}\}$ may be \emph{renormalized}, or rescaled
to yield a new
dynamical system of the same
general shape as that of the original map $r(z)$?
\par
Since the scaling equation (\ref{eqscaling0}) from wavelet theory is a
renormalization, our application of multiresolutions to the Julia sets $X(r)$
of complex dynamics suggests a useful approach to renormalization in this
context. The drawback of this approach is that it relies on a rather
unwieldy Hilbert space built on $X(r)$, or rather on the projective system
$X(r)_\infty$ built in turn over $X(r)$ . So we are left with translating our
Hil\-bert-space-the\-o\-ret\-ic normalization back to the direct and geometric
algorithms on the complex plane.
\par
     General and rigorous results from complex dynamics state that under
additional geometric hypotheses, renormalized dynamical systems range in a
compact family.

     The use of geometric tools from Hilbert space seems promising for renormalization questions (see, e.g.,
\cite{MR2001c:60064}, \cite{Blu04}, \cite{BreJo}, and \cite{VoYu})
since notions of repetition up to similarity of form at infinitely
many scales are common in the Hilbert-space approach to
multiresolutions; see, e.g., \cite{Don04}. And this
self-similarity up to scale parallels a basic feature of
renormalization questions in physics: for a variety of instances
of dynamics in physics, and in telecommunication
\cite{MR2001c:60064,Blu04,BJKR01}, we encounter scaling laws of
self-similarity; i.e., we observe that a phenomenon reproduces
itself on different time and/or space scales.

Self-similar processes are stochastic processes that are invariant in
distribution under suitable scaling of time and/or space (details below!)
Fractional Brownian motion is perhaps the best known of these, and it is used
in telecommunication and in stochastic integration. While the underlying idea
behind this can be traced back to Kolmogorov, it is only recently, with the
advent of wavelet methods, that its computational power has come more into
focus, see e.g., \cite{MR2001c:60064}. But at the same time, this connection to wavelet
analysis is now bringing the \emph{operator-theoretic features} in the subject to
the fore.

    In statistics, we observe that fractional Brownian motion (a Gaussian
process $B(t)$ with $E\{B(t)\} = 0$, and covariance $E\{B(t)B(s)\}$ given by a certain $h$-fractional law) has the property that there is a number $h$ such that, for all $a$, the two processes $B(at)$ and $a^h B(t)$ have the same finite-dimensional distributions. The scaling feature of fractional Brownian motion, and of its corresponding white-noise process, is used in telecommunication, and in stochastic integration; and this connection has been one of the new and more exciting domains of applications of wavelet tools%
, both pure and applied (Donoho, Daubechies, Meyer, etc.
\cite{Don04,MR2004a:42046,MR2003d:46102}).

This is a new direction of pure mathematics which
makes essential contact with problems that are not normally investigated
in the context of harmonic analysis.

In our proofs, we take advantage of our Hilbert-%
space formulation
of the notion of
\label{sutsa}\emph{similarity up to scale}, and the use of scales of closed subspaces
$\mathcal{V}_n$ in a Hilbert space $\mathcal{H}$.
The unitary operator $U$ which scales between the spaces may arise from a
substitution in nonlinear dynamics, such as in the context of
Julia sets \cite{Be91}; it may ``scale'' between  the sigma-algebras in
a martingale \cite{Doo49,Doo61,Wil02}; it may be a scaling of large volumes of
data (see, e.g., \cite{AFTV00}); it may be the scaling operation in a
 fractional Brownian motion \cite{MR2001c:60064,Lin93}; it may be cell averages
in finite elements \cite{MR2004c:65178}; or it may
be the traditional dyadic scaling operator in $\mathcal{H}=L^2(\mathbb{R})$
in wavelet frame analysis \cite{CR95,MR2004c:65178}.

Once the system\label{Otsb} $(\mathcal{H}, (\mathcal{V}_n), U)$ is given, then it  follows that there is a
unitary representation $V$ of $\mathbb{T}$ in $\mathcal{H}$  which defines a grading of operators
on
$\mathcal{H}$. Moreover $U$ then defines a similarity up to scale precisely when $U$ has
grade one, see 
(\ref{eqPJ820.1}) for definition. We also  refer to the papers
\cite{Arv02b,BaVi02} for details and applications of this  notion in operator theory.

The general idea of assigning degrees to operators in Hilbert  space has
served as a powerful tool in other areas of mathematical  physics (see,
e.g.,
\cite{Ar97}, \cite{BrRoII}; and \cite{BoCo95} on phase  transition problems); operator
theory \cite{Arv02a,Arv02b,MR2004i:46109,MR2004j:46034}; and operator algebras; see, e.g.,  \cite{Con94,CoSk02,BCE95}.

Our research, e.g., \cite{DuJo04a}, in fact already indicates that the
nonlinear problems sketched above can be attacked with the use of our
operator-theoretic framework.

\subsection{\label{crossed}Wavelet representations of operator-algebraic crossed products}

The construction of wavelets requires a choice of a low-pass filter $m_0\in L^\infty(\mathbb{T})$ and of $N-1$ high-pass filters $m_i\in L^\infty(\mathbb{T})$ that satisfy certain orthogonality conditions. As shown in \cite{BrJo97b}, these choices can be put in one-to-one correspondence with a class of covariant representations of the Cuntz algebra $\mathcal{O}_N$. In the case of an arbitrary $N$-to-one dynamical system $r \colon X\rightarrow X$, the wavelet construction will again involve a careful choice of the filters, and a new class of representations is obtain. The question here is to classify these representations and see how they depend on the dynamical system $r$. A canonical choice of the filters could provide an invariant for the dynamical system.
\par
In a recent paper \cite{Exe03}, the crossed-product of a $C^*$-algebra by an endomorphism is constructed using the transfer operator. The required covariant relations are satisfied by the operators that we introduced in \cite{DuJo04a}, \cite{DuJo04b}.
\par
There are two spaces for these representations: one is the core space of the multiresolution $V_0=L^2(X,\mu)$. Here we have the abelian algebra of multiplication operators and the isometry $S \colon f\mapsto m_0\,f\circ r$. Together they provide representations of the crossed product by an endomorphism.
The second space is associated to a dilated measure $\hat\mu$ on the solenoid of $r$, $L^2(X_\infty,\hat\mu)$ (see \cite{DuJo04a}, \cite{DuJo04b}). The isometry $S$ is dilated to a unitary $\hat S$ which preserves the covariance. In this case, we are dealing with representations of crossed products by automorphisms.
\par
Thus we have four objects interacting with one another: crossed products by endomorphisms and their representations, and crossed products by automorphisms and its representations on $L^2(\hat X,\hat\mu)$.
\par
The representations come with a rich structure given by the multiresolution, scaling function and wavelets. Therefore their analysis and classification seems promising.

\par

\subsection{\label{KMS}KMS-states and equilibrium measures}
In \cite{Exe04}, the KMS states on the Cuntz--Pimsner algebra for a certain one-parameter group of automorphisms are obtained by composing the unique equilibrium measure on an abelian subalgebra (i.e., the measure which is invariant to the transfer operator) with the conditional expectation. These results were further generalized in \cite{KuRe04} for expansive maps and the associated groupoid algebras.
\par
As we have seen in \cite{DuJo03} for the case when $r(z)=z^N$, the covariant representations associated to some low-pass filter $m_0$ are highly dependent on the equilibrium measure $\nu$. 
As outlined above, we consider here an expansive dynamical system built on a
given finite-to-one mapping $r\colon X\rightarrow X$, and generalizing the familiar case
of the winding mapping $z\rightarrow z^N$ on $\mathbb{T}$. In our $r\colon X\rightarrow  X$  context, we then
study weight functions $V \colon X\rightarrow  [0, \infty)$. If the zero-set of $V$ is now
assumed finite, we have so far generalized what is known in the $z\rightarrow z^N$
case. Our results take the form of a certain dichotomy for the admissible
Perron--Frobenius--Ruelle measures, and they fit within the framework of \cite{KuRe04}.
Hence our construction also yields covariant representations and symmetry in
the context of \cite{KuRe04}. Our equilibrium measures in turn induce the kind of
KMS-states studied in \cite{KuRe04}.
\par
The results of \cite{KuRe04} and \cite{Exe04} restrict the weight function (in our case represented by the absolute square of the low-pass filter, $|m_0|^2$) to being strictly positive and H\"older continuous. These restrictions are required because of the form of the Ruelle--Perron--Frobenius theorem used, which guarantees the existence and uniqueness of the equilibrium measure.
\par
However, the classical wavelet theory shows that many interesting examples require a ``low-pass condition'', $m_0(1)=\sqrt{N}$, which brings zeroes for $|m_0|^2$.
\par
Our martingale approach is much less restrictive than the conditions of \cite{KuRe04}: it allows zeroes and discontinuities. Coupled with a recent, more general form of Ruelle's theorem from \cite{FaJi01a}, we hope to extend the results of \cite{KuRe04} and be able to find the KMS states in a more general case. Since the existence of zeroes can imply a multiplicity of equilibrium measures (see \cite{Dut04a}, \cite{BrJo02b}) a new phenomenon might occur such as spontaneously breaking symmetry.

\par
For each KMS state, one can construct the GNS representation. Interesting results are known \cite{Oka03} for the subshifts of finite type and the Cuntz--Krieger algebra, when the construction yields type $\mathrm{III}_\lambda$ AFD factors.
\par
The natural question is what type of factors appear in the general case of $r \colon X\rightarrow X$. Again, a canonical choice for the low-pass filter can provide invariants for the dynamical system.
\par
Historically, von-Neumann-algebra techniques have served
as powerful tools for analyzing discrete structures. Our work so far suggests
that the iteration systems $(X,r)$ will indeed induce interesting von-Neumann-%
algebra isomorphism classes.

\subsection{\label{Dim}Dimension groups} We propose to use our geometric
algorithms \cite{BJKR01,BJKR02} for deciding order isomorphism of
dimension groups for the purpose of deciding isomorphism questions
which arise in our multiresolution analysis built on Julia sets.
The project is to use the results in
\cite{BJKR01} 
to classify substitution tilings
in the sense of \cite{Rad91}, \cite{Rad99}.
\par

\subsection{\label{Riesz}Equilibrium measures, harmonic functions for the transfer
operator, and infinite Riesz products}
In \cite{BrJo02b}, \cite{Dut04a}, for the dynamical system $z\mapsto z^N$, we were able to compute the equilibrium measures, the fixed points of the transfer operator and its ergodic decomposition by analyzing the commutant of the covariant representations associated to the filter $m_0$. In \cite{DuJo04a}, \cite{DuJo04b} we have shown that also in the general case there is a one-to-one correspondence between harmonic functions for the transfer operator and the commutant. Thus an explicit form of the covariant representation can give explicit solutions for the eigenvalue problem $R_{m_0}h=h$, which are of interest in ergodic theory (see \cite{Bal00}).
\par
In some cases, the equilibrium measures are Riesz products (see \cite{DuJo03}); these are examples of exotic measures which arise in harmonic analysis \cite{BBK04,Erd40,Kah71,MR2002k:93007,Rie18,Bis95,Bro75,BrMo75,Kat87,Mey79,Rit79}. The wavelet-operator-theoretic approach may provide new insight into their properties.

\subsection{\label{Ruelle}Non-uniqueness of the Ruelle--Perron--Frobenius data}
A substantial part of the current literature in the Ruelle--Perron--Frobenius operator in its many guises is primarily about conditions for uniqueness of KMS; so that means no phase transition. This is ironic since Ruelle's pioneering work was motivated by phase-transition questions, i.e., non-uniqueness.

In fact non-uniqueness is much more interesting in a variety of applications: That is when we must study weight functions $W$ in the Ruelle operator $R = R_W$ when the standard rather restrictive conditions on $W$ are in force. The much celebrated Ruelle theorem for the transition operator $R$ gives existence and uniqueness of Perron--Frobenius data under these restrictive conditions. But both demands from physics, and from wavelets, suggest strong and promising research interest in relaxing the stringent restrictions that are typical in the literature, referring to assumptions on $W$.  There is a variety of exciting possibilities.  They would give us existence, but not necessarily uniqueness in the conclusion of a new Ruelle-type theorem.

\subsection{\label{Indu}Induced measures}
A basic tool in stochastic processes (from probability theory)
involves a construction on a ``large'' projective space $X_\infty$,
based on some marginal measure on some coordinate space $X$. Our projective
limit space $X_\infty$ will be constructed from a finite
branching process.
\par
\par
Let $A$ be a $k\times k$ matrix with entries in $\{0,1\}$.
Suppose every column in $A$ contains an entry $1$.
\par
Set
$X(A):=\left\{\,(\xi_i)_{i\in\bn}\in\prod_{\bn}\{1,\dots,k\}\bigm|A(\xi_i,\xi_{i+1})=1\,\right\}$
and
$$r_A(\xi_1,\xi_2,\dots)=(\xi_2,\xi_3,\dots)\mbox{\quad for }\xi\in X(A).$$
Then $r_A$ is a subshift, and the pair $(X(A),r_A)$ satisfies our
conditions.
\par
It is known \cite{Rue89} that, for each $A$, as described, the
corresponding system $r_A \colon X(A)\rightarrow X(A)$ has a unique
strongly $r_A$-invariant probability measure, $\rho_A$, i.e., a
probability measure on $X(A)$ such that $\rho_A\circ R_1=\rho_A$, where $R_1$ is defined as in (\ref{eq2_18}) below.
\par
We now turn to the connection between measures on $X$ and an associated family of induced measures on $X_\infty$, and we characterize those
measures $X_\infty$ which are quasi-invariant with respect to the
invertible mapping $\hat r$,
where
$
X_\infty:=\{\,\hat
x=(x_0,x_1,\dots)\in\prod_{\bn_0}X\mid r(x_{n+1})=x_n\,\}
$,
$\hat r(\hat x)=(r(x_0),x_0,x_1,\dots)$.

\par
In our extension of measures from $X$ to $X_\infty$, we must keep
track of the transfer from one step to the next, and there is an
operator which accomplishes this, Ruelle's transfer operator
\begin{equation}
R_Wf(x)=\frac{1}{\#\,r^{-1}(x)}\sum_{r(y)=x}W(y)f(y),\qquad f\in
L^1(X,\mu).
\label{eq2_18}
\end{equation}
In its original form it was introduced in
\cite{Rue89}, but since that, it has found a variety of
applications, see e.g., \cite{Bal00}. For use of the Ruelle
operator in wavelet theory, we refer to \cite{Jor01a} and
\cite{Dut04a}.

\par
The Hilbert spaces of functions on
$X_\infty$ will be realized as a Hilbert spaces of martingales.
This is consistent with our treatment of wavelet resolutions as
martingales. This was first suggested in \cite{Gun99} in
connection with wavelet analysis.
\par
 In
\cite{DuJo04a}, we studied the following restrictive setup: we
assumed that $X$ carries a probability measure $\mu$ which is {\it
strongly $r$-invariant}. By this we mean that
\begin{equation}\label{eq2_17}
\int_Xf\,d\mu=\int_X\frac{1}{\#\,r^{-1}(x)}\sum_{y\in X,
r(y)=x}f(y)\,d\mu(x),\qquad f\in L^\infty(X).
\end{equation}
\par
If, for example $X=\Br/\bz$, and $r(x)=2x\bmod 1$, then the Haar
measure on $\Br/\bz={}$Lebesgue measure on $[0,1)$ is the unique
strongly $r$-invariant measure on $X$.
\par
Suppose the weight function $V$ is bounded and measurable. Then
define $R=R_V$, the {\it Ruelle operator}, by formula
(\ref{eq2_18}), with $V=W$.

\begin{theorem}\label{th4_18}\textup{(\cite{DuJo04a})} Let $r \colon X\rightarrow X$ and $\xir$
be as described above, and suppose that $X$ has a strongly
$r$-invariant measure $\mu$. Let $V$ be a non-negative, measurable
function on $X$, and let $R_V$ be the corresponding Ruelle
operator.
\begin{enumerate}
\item There is a unique measure $\hat\mu$ on $\xir$ such that
\begin{enumerate}
\item $\hat\mu\circ\theta_0^{-1}\ll\mu$ \textup{(}set
$h=\frac{d(\hat\mu\circ\theta_0^{-1})}{d\mu}$\textup{)}, \item
$\int_Xf\,d\hat\mu\circ\theta_n^{-1}=\int_XR^n_V(fh)\,d\mu,\qquad n\in\bn_0.$
\end{enumerate}
\item The measure $\hat\mu$ on $\xir$ satisfies
\begin{equation}\label{eq2_19}
\frac{d(\hat\mu\circ\hat r)}{d\hat\mu}=V\circ\theta_0,
\end{equation}
and
\begin{equation}\label{eq2_20}
R_Vh=h.
\end{equation}
\end{enumerate}
\end{theorem}

\par
If the function $V \colon X\rightarrow[0,\infty)$ is given, we define
$$V^{(n)}(x):=V(x)V(r(x))\cdots V(r^{n-1}(x)),$$
and set $d\mu_n:=V^{(n)}d\mu_0$. Our result states that the
corresponding measure $\hat\mu$ on $\xir$ is $V$-quasi-invariant
if and only if
\begin{equation}\label{eq5_no}
d\mu_0=(V\,d\mu_0)\circ r^{-1}.
\end{equation}

\begin{theorem}\label{th4_19}
Let $V \colon X\rightarrow[0,\infty)$ be $\mathfrak{B}$-measurable, and
let $\mu_0$ be a measure on $X$ satisfying the following
fixed-point property
\begin{equation}\label{eq5_1}
d\mu_0=(V\,d\mu_0)\circ r^{-1}.
\end{equation}
Then there exists a unique measure $\hat\mu$ on $\xir$ such that
\begin{equation}\label{eq5_2}
\frac{d(\hat\mu\circ\hat r)}{d\hat\mu}=V\circ\theta_0
\end{equation}
and
$$\hat\mu\circ\theta_0^{-1}=\mu_0.$$
\end{theorem}

\par
\begin{definition}\label{def6_no}
Let $V \colon X\rightarrow[0,\infty)$ be bounded and
$\mathfrak{B}$-measurable. We use the notation
$$M^V(X):=\{\,\mu\in M(X)\mid d\mu=(V\,d\mu)\circ r^{-1}\,\}.$$
\par
For measures $\hat\mu$ on $(\xir,\mathfrak{B}_\infty)$ we
introduce
\begin{equation}\label{eqmqi}M_{qi}^V(\xir):=\left\{\,\hat\mu\in M(\xir)\Bigm|\hat\mu\circ\hat
r\ll\hat\mu\mbox{ and }\frac{d(\hat\mu\circ\hat
r)}{d\hat\mu}=V\circ\theta_0\,\right\}.\end{equation}
\end{definition}
As in Definition \ref{def6_no}, let $X$, $r$, and $V$ be as before, i.e., $r$ is a
finite-to-one endomorphism of $X$, and $V  \colon  X \rightarrow [0, \infty)$ is a given
weight function. In the next theorem, we establish a canonical bijection
between the convex set $M^V(X)$ of measures on $X$ with the set of
$V$-quasi-invariant measures on $X_\infty$, which we call
$M_{qi}^V(X_\infty(r))$, see (\ref{eqmqi}).

\par
The results so far in this section may be summarized as follows:
\begin{theorem}\label{th6_no_no}
Let $V$ be as in Definition \ref{def6_no}. For measures $\hat\mu$
on $\xir$ and $n\in\bn_0$, define
$$C_n(\hat\mu):=\hat\mu\circ\theta_n^{-1}.$$
Then $C_0$ is a bijective affine isomorphism of $M_{qi}^V(\xir)$
onto $M^V(X)$ that preserves the total measure, i.e.,
$C_0(\hat\mu)(X)=\hat\mu(\xir)$ for all $\hat\mu\in
M_{qi}^V(\xir)$.
\end{theorem}

\begin{theorem}\label{th6_nooo}
Let $V \colon X\rightarrow[0,\infty)$ be continuous. Assume that there
exist some measure $\nu$ on $(X,\mathfrak{B})$ and two numbers
$0<a<b$ such that
\begin{equation}\label{eq6_nooo_1}
a\leq\nu(X)\leq b, \mbox{ and }a\leq\int_X V^{(n)}\,d\nu\leq
b\mbox{ for all }n\in\bn.
\end{equation}
Then there exists a measure $\mu_0$ on $(X,\mathfrak{B})$ that
satisfies
$$d\mu_0=(V\,d\mu_0)\circ r^{-1},$$
and there exists a $V$-quasi-invariant measure $\hat\mu$ on
$(\xir,\mathfrak{B}_\infty)$.
\end{theorem}

\begin{theorem}\label{th6_7}
Let $(X,\mathfrak{B})$, and $r \colon X\rightarrow X$, be as described
above. Suppose $V \colon X\rightarrow[0,\infty)$ is measurable,
$$\frac{1}{\#\,r^{-1}(x)}\sum_{r(y)=x}V(y)\leq1,$$
and that some probability measure $\nu_V$ on $X$ satisfies
\begin{equation}\label{eq6_7_no}
\nu_V\circ R_V=\nu_V.
\end{equation} Assume also that
$(X,\mathfrak{B})$ carries a strongly $r$-invariant probability
measure $\rho$, such that
\begin{equation}\label{eq6_7_2}
\rho(\{\,x\in X\mid V(x)>0\,\})>0.
\end{equation}
Then
\begin{enumerate}
\item $T_V^n(\,d\rho)=R_V^n(\mathbf{1})\,d\rho,$ for $n\in\bn,$
where $\mathbf{1}$ denotes the constant function one. \item
 \textup{[}Monotonicity\/\textup{]} $\quad \cdots\leq R_V^{n+1}(\mathbf{1})\leq
R_V^n(\mathbf{1})\leq\dots\leq\mathbf{1}$, pointwise on $X$.

\item The limit $\lim_{n\rightarrow}R_V^n(\mathbf{1})=h_V$
exists, $R_Vh_V=h_V$, and
\begin{equation}\label{eq6_7_3}
\rho(\{\,x\in X\mid h_V(x)>0\,\})>0.
\end{equation}
\item The measure $d\mu_0^{(V)}=h_V\,d\rho$ is a solution to the
fixed-point problem $$T_V(\mu_0^{(V)})=\mu_0^{(V)}.$$ \item The
sequence $d\mu_n^{(V)}=V^{(n)}\,h_V\,d\rho$ defines a unique
$\hat\mu^{(V)}$ as in Theorem \ref{th4_18}; and \item
$\mu_n^{(V)}(f)=\int_XR_V^n(fh_V)\,d\rho$ for all bounded
measurable functions $f$ on $X$, and all $n\in\bn$.
\end{enumerate}
\setcounter{saveenumi}{\value{enumi}}
Finally, 
\begin{enumerate}
\setcounter{enumi}{\value{saveenumi}}
\item The measure $\hat\mu^{(V)}$ on $\xir$ satisfying
$\hat\mu^{(V)}\circ\theta_n^{-1}=\mu_n^{(V)}$ has total mass
$$\hat\mu^{(V)}(\xir)=\rho(h_V)=\int_Xh_V(x)\,d\rho(x).$$
\end{enumerate}
\end{theorem}

\begin{definition}\label{def7_1_no}
A measure $\mu_0\in M_1^V(X)$ is called {\it relatively ergodic}
with respect to $(r,V)$ if the only non-negative, bounded
$\mathfrak{B}$-measurable functions $f$ on $X$ satisfying
$$E_{\mu_0}(Vf)=E_{\mu_0}(V)f\circ r,\,\mbox{ pointwise }\mu_0\circ
r^{-1}\mbox{-a.e.},$$ are the functions which are constant
$\mu_0$-a.e.
\end{definition}

\par

Since we have a canonical bijection between the two compact convex sets
of measures $M_1^V(X)$ and $M_{qi,1}^V(X_\infty(r))$, the natural question
arises as to the extreme points. This is answered in our next theorem. We show that there are notions ergodicity for each of the two
spaces $X$ and $X_\infty(r)$ which are equivalent to extremality in the
corresponding compact convex sets of measures.

\begin{theorem}\label{th7_2}
Let $V \colon X\rightarrow[0,\infty)$ be bounded and measurable. Let
$$\hat\mu\in M_{qi,1}^V(\xir)\text{, and
}\mu_0:=\hat\mu\circ\theta_0^{-1}\in M_1^V(X).$$ The following
assertions are equivalent:
\begin{enumerate}
\item $\hat\mu$ is an extreme point of $M_{qi,1}^V(\xir)$; \item
$V\circ\theta_0\,d\hat\mu$ is ergodic with respect to $\hat r$;
\item $\mu_0$ is an extreme point of $M_1^V(X)$; \item $\mu_0$ is
relatively ergodic with respect to $(r,V)$.
\end{enumerate}
\end{theorem}

We now turn to the next two theorems.
These are counterparts of our di\-men\-sion-coun\-ting functions which we
outlined above in connection with Theorem \ref{thmul1}; see especially Remark \ref{rem4_9}.
They concern the correspondence
between the two classes of measures, certain on $X$ (see Theorem \ref{th7_2}), and the
associated induced measures on the projective space $X_\infty$. Our proofs
are technical and will not be included. (The reader is referred to \cite{DuJo04c} for
further details.)  Rather we only give a few suggestive hints: Below we
outline certain combinatorial concepts and functions which are central for
the arguments. Since they involve a counting principle, they have an
intuitive flavor which we hope will offer some insight into our theorems.

\par
Let $X$ be a non-empty set, $\mathfrak{B}$ a sigma-algebra of
subsets of $X$, and $r \colon X\rightarrow X$ an onto, finite-to-one, and
$\mathfrak{B}$-measurable map.
\par
We will assume in addition that we can label measurably the
branches of the inverse of $r$. By this, we mean that the
following conditions are satisfied:
\begin{equation}\label{eq2_1}
\mbox{The map}\quad \mathfrak{c} \colon X\ni
x\mapsto\#\,r^{-1}(x)<\infty\mbox{ is measurable }.
\end{equation}
We denote by $A_i$ the set
$$A_i:=\{\,x\in X\mid \mathfrak{c}(x)=\#\,r^{-1}(x)\geq
i\,\},\qquad i\in\bn.$$ Equation (\ref{eq2_1}) implies that the
sets $A_i$ are measurable. Also they form a decreasing sequence
and, since the map is finite-to-one,
$$X=\bigcup_{i=1}^\infty(A_{i+1}\setminus A_i).$$
Then, we assume that there exist measurable maps
$\tau_i \colon A_i\rightarrow X$, $i\in\{1,2,\dots\}$ such that
\begin{equation}\label{eq2_2}
r^{-1}(x)=\{\tau_1(x),\dots,\tau_{\mathfrak{c}(x)}(x)\},\qquad
x\in X,
\end{equation}
\begin{equation}\label{eq2_3}
\tau_i(A_i)\in\mathfrak{B}\mbox{\quad for all }i\in\{1,2,\dots\}.
\end{equation}
Thus $\tau_1(x),\dots,\tau_{\mathfrak{c}(x)}(x)$ is a list without
repetitions of the ``roots'' of $x$, $r^{-1}(x)$.
\par
{}From (\ref{eq2_2}) we obtain also that
\begin{equation}\label{eq2_4}
\tau_i(A_i)\cap\tau_j(A_j)=\emptyset,\mbox{ if }i\neq j,
\end{equation}
and
\begin{equation}\label{eq2_5dup}
\bigcup_{i=1}^\infty\tau_i(A_i)=X.
\end{equation}
\par In the sequel, we will use the following notation: for a
function $f \colon X\rightarrow\bc$, we denote by $f\circ\tau_i$ the function
$$f\circ\tau_i(x):=\left\{\begin{array}{ll}
f(\tau_i(x))&\mbox{ if }x\in A_i,\\
0&\mbox{ if }x\in X\setminus A_i.\end{array}\right.
$$

Our Theorem \ref{thpr_1} depends on the existence of some strongly invariant
measure $\mu$ on $X$, when the system $(X, r)$ is given. However, in the general
measurable category, such a strongly invariant measure $\mu$ on $X$ may in
fact not exist; or if it does, it may not be available by computation. In
fact, recent wavelet results (for frequency localized wavelets, see, e.g.,
\cite{BJMP03b} and \cite{BJMP04}) suggest the need for theorems in the more general class of
dynamical systems $(X, r)$.
\par

In the next theorem (Theorem \ref{th4_14}), we provide for each system, $X$, $r$, and $V$ a
substitute for the existence of strongly invariant measures. We show that
there is a certain fixed-point property for measures on $X$ which depends
on $V$ but not on the a priori existence of strongly $r$-invariant measures,
and which instead is determined by a certain modular function $\Delta$ on $X$.
This modular function in turn allows us to prove a dis-integration
theorem for our $V$-quasi-invariant measures $\hat\mu$ on $X_\infty$. In Theorem
\ref{th4_17}, we give a formula for this dis-integration of a $V$-quasi-invariant
measure $\hat\mu$ on $X_\infty$ in the presence of a modular function $\Delta$.
Our dis-integration of $\hat\mu$ is over a Markov process $P_x$, for $x$ in $X$,
which depends on $\Delta$, but otherwise is analogous to the process $P_x$ we
used in Proposition \ref{prop3_1}.

\begin{theorem}\label{th4_14}
Let $(X,\mathfrak{B})$ and $r \colon X\rightarrow X$ be as above. Let
$V \colon X\rightarrow [0,\infty)$ be a bounded
$\mathfrak{B}$-measurable map. For a measure $\mu_0$ on
$(X,\mathfrak{B})$, the following assertions are equivalent.
\begin{enumerate}
\item\label{th4_14(1)} The measure $\mu_0$ has the fixed-point property
\begin{equation}\label{eq2_1_1}
\int_XV\,f\circ r\,d\mu_0=\int_Xf\,d\mu_0,\mbox{ for all }f\in
L^\infty(X,\mu_0).
\end{equation}
\item\label{th4_14(2)} There exists a non-negative, $\mathfrak{B}$-measurable map
$\Delta$ \textup{(}depending on $V$ and $\mu_0$\textup{)} such that
\begin{equation}\label{eq2_1_2}
\sum_{r(y)=x}\Delta(y)=1,\mbox{ for }\mu_0\mbox{-a.e. }x\in X,
\end{equation}
and
\begin{equation}\label{eq2_1_3}
\int_XV\,f\,d\mu_0=\int_X\sum_{r(y)=x}\Delta(y)f(y)\,d\mu_0(x),\mbox{
for all }f\in L^\infty(X,\mu_0).
\end{equation}
\end{enumerate}
Moreover, when the assertions are true, $\Delta$ is unique up
to $V\,d\mu_0$-measure zero.
\end{theorem}

\par

We recall the definitions: if $\mathfrak{B}$ is a sigma-algebra
on a set $X$ and $r \colon X\rightarrow X$ is a finite-to-one, onto and
measurable map, we denote by $X_\infty$ the set
$$\xir:=\left\{\,(x_0,x_1,\dots)\in\prod_{n\in\bn_0}X\biggm|r(x_{n+1})=x_n,\mbox{
for all }n\in\bn_0\,\right\}.$$ We denote the projections by
$\theta_n \colon \xir\ni(x_0,x_1,\dots)\mapsto x_n\in X.$ The union of the
pull-backs $\theta_n^{-1}(\mathfrak{B})$ generates a
sigma-algebra $\mathfrak{B}_\infty$. The map $r$ extends to a
measurable bijection $\hat r \colon \xir\rightarrow\xir$ defined by
$$\hat r(x_0,x_1,\dots)=(r(x_0),x_0,x_1,\dots).$$
\par
Let $V \colon X\rightarrow[0,\infty)$ be a measurable, bounded function.
We say that a measure $\mu_0$ on $(X,\mathfrak{B})$ has the {\it
fixed-point property} if
$$\int_XV\,f\circ r\,d\mu_0=\int_Xf\,d\mu_0,\qquad f\in
L^\infty(X,\mu_0).$$ We say that a measure $\hat\mu$ on
$(\xir,\mathfrak{B}_\infty)$ is {\it $V$-quasi-invariant} if
$$d(\hat\mu\circ\hat r)=V\circ\theta_0\,d\hat\mu.$$
We recall the following result from \cite{DuJo04b}.
\begin{theorem}\label{th4_15}
There exists a one-to-one correspondence between measures $\mu_0$
on $X$ with the fixed-point property and $V$-quasi-invariant
measures $\hat\mu$ on $\xir$, given by
$$\mu_0=\hat\mu\circ\theta_0^{-1}.$$
\end{theorem}

\par
\begin{proposition}\label{prop4_16}
Let $(X,\mathfrak{B})$, $r \colon X\rightarrow X$ and be as above, and
let $D \colon X\rightarrow[0,\infty)$ be a measurable function with the
property that
\begin{equation}\label{eq2_3_1}
\sum_{r(y)=x}D(y)=1.
\end{equation}
Denote by $D^{(n)}$ the product of compositions
\begin{equation}\label{eq2_3_2}
D^{(n)}:=D\cdot D\circ r\cdot\,\cdots\,\cdot D\circ
r^{n-1},\qquad n\in\bn,\quad D^{(0)}:=1.
\end{equation} Then for each $x_0\in X$, there exists a Radon probability measure $P_{x_0}$
on $\Omega_{x_0}$ such that, if $f$ is a bounded measurable
function on $\Omega_{x_0}$ which depends only on the first $n$
coordinates $x_1,\dots,x_n$, then
\begin{equation}\label{eq2_3_3}
\int_{\Omega_{x_0}}
f(\omega)\,dP_{x_0}(\omega)=\sum_{r^n(x_n)=x_0}D^{(n)}(x_n)f(x_1,\dots,x_n).
\end{equation}
\end{proposition}

\par
\begin{theorem}\label{th4_17}
Let $(X,\mathfrak{B})$, $r \colon X\rightarrow X$ and
$V \colon X\rightarrow[0,\infty)$ be as above. Let $\mu_0$ be a measure
on $(X,\mathfrak{B})$ with the fixed-point property
\textup{(\ref{eq2_1_1})}. Let $\Delta \colon X\rightarrow[0,1]$ be the function
associated to $V$ and $\mu_0$ as in Theorem \ref{th4_14}, and let
$\hat\mu$ be the unique $V$-quasi-invariant measure on $\xir$
with
$$\hat\mu\circ\theta_0^{-1}=\mu_0,$$
as in Theorem \ref{th4_15}. For $\Delta$ define the measures
$P_{x_0}$ as in Proposition \ref{prop4_16}. Then, for all bounded
measurable functions $f$ on $\xir$,
\begin{equation}\label{eq2_4_1}
\int_{\xir}f\,d\hat\mu=\int_X\int_{\Omega _{x_0}}
f(x_0,\omega)\,dP_{x_0}(\omega)\,d\mu_0(x_0).
\end{equation}
\end{theorem}

\subsection{\label{Haus}Hausdorff measure and wavelet bases}
In this problem we propose wavelet bases in the context of Hausdorff
measure of fractional dimension between 0 and 1. While our fractal wavelet
theory has points of similarity that it shares with the standard case of
Lebesgue measure on the line, there are also sharp contrasts.

It is well known that the Hilbert spaces $L^{2}(\mathbb{R})$ has a rich
family of orthonormal bases of the following form:
$
\psi _{j,k}(x)=2^{j/2}\psi (2^{j}x-k),\, j,k\in \mathbb{Z},
$
where $\psi $ is a single function $\in L^{2}(\mathbb{R})$.
Take for example Haar's function
$\psi (x)=\chi _{I}(2x)-\chi _{I}(2x-1)$
where $I=[0,1]$ is the unit interval.
Let $\mathbf{C}$ be the middle-third Cantor set. Then the corresponding
indicator function $\varphi _{\mathbf{C}}:=\chi _{\mathbf{C}}$
satisfies the scaling identity (see (\ref{eqscaling0})),
$
\varphi _{\mathbf{C}}(\frac{x}{3})=\varphi _{\mathbf{C}}(x)+\varphi _{%
\mathbf{C}}(x-2)\text{.}$

In \cite{DuJo03} we use this as the starting point for a
multiresolution construction in a separable Hilbert space
built with Hausdorff measure, and we identify the two mother functions
which define the associated wavelet ONB.

Since both constructions, the first one for the Lebesgue measure, and the
second one for the Hausdorff version $(dx)^{s}$, arise from scaling and
subdivision, it seems reasonable to expect multiresolution wavelets also in
Hilbert spaces constructed on the scaled Hausdorff measures $\mathcal{H}^{s}$
which are basic for the kind of iterated function systems which give Cantor
constructions built on scaling and translations by lattices.

\begin{acknowledgements}
The authors gratefully acknowledge constructive discussions 
with Richard Gundy and Roger Nussbaum. Further, we are pleased to thank Brian 
Treadway for help with \TeX, with style files, and with proofreading. And 
especially for making a number of very helpful suggestions. Both authors 
have benefited from discussions with the members of our NSF Focused Research 
Group (DMS-0139473), ``Wavelets, Frames, and Operator Theory,'' especially 
from discussions with Larry Baggett, John Benedetto, David Larson, and Gestur \'Olaffson.  
The second named author is grateful for helpful discussions and suggestions from 
Kathy Merrill, Wayne Lawton, Paul Muhly, Judy Packer, Iain Raeburn, and Wai-Shing Tang. 
Several of these were participants and organizers of an 
August 4--7, 2004, Wavelet Workshop 
(``Functional and harmonic analysis of wavelets and frames'') 
at The National University of Singapore (NUS). And we benefited from inspiration and
enlightening discussions with this workshop as well. The second named author is grateful 
to the hosts at NUS in Singapore for support and generous hospitality in August of 2004.
Finally, we thank Daniel Alpay for inviting us to write this paper, and for encouragements.
\end{acknowledgements}

\providecommand{\bysame}{\leavevmode\hbox
to3em{\hrulefill}\thinspace}
\providecommand{\MR}{\relax\ifhmode\unskip\space\fi MR }
\renewcommand{\MR}[1]{\relax}
\providecommand{\MRhref}[2]{%
  \href{http://www.ams.org/mathscinet-getitem?mr=#1}{#2}
} \providecommand{\href}[2]{#2}



\end{document}